\documentclass[11pt]{amsart}
\usepackage[margin=1.1in]{geometry}
\usepackage{amsmath,amssymb,amsthm,mathtools}
\usepackage{needspace}
\usepackage{cite}
\usepackage{xcolor}
\definecolor{linkink}{RGB}{0,0,0}
\usepackage[colorlinks=true,allcolors=linkink,pdftitle={Weak-strong uniqueness for the Landau equation}]{hyperref}
\newtheorem{theorem}{Theorem}[section]
\newtheorem{proposition}[theorem]{Proposition}
\newtheorem{lemma}[theorem]{Lemma}
\newtheorem{corollary}[theorem]{Corollary}
\theoremstyle{definition}
\newtheorem*{definition*}{Definition}
\theoremstyle{remark}
\newtheorem{remark}[theorem]{Remark}
\numberwithin{equation}{section}
\newcommand{\R}{\mathbb{R}}
\newcommand{\ba}{\bar a}
\newcommand{\bc}{\bar c}
\newcommand{\jv}{\langle v\rangle}
\newcommand{\jy}{\langle y\rangle}
\newcommand{\K}{\mathcal{K}}
\newcommand{\cP}{\mathcal{P}}
\newcommand{\cZ}{\mathcal{Z}}
\newcommand{\tr}{\operatorname{tr}}
\newcommand{\dd}{\,\mathrm{d}}
\newcommand{\tv}{\tilde v}
\newcommand{\tx}{\tilde x}
\newcommand{\ttil}{\tilde t}
\newcommand{\te}{\theta_\epsilon}
\newcommand{\step}[1]{\smallskip\noindent\textbf{#1}\ }
\DeclareRobustCommand{\rev}[1]{\leavevmode{\color{black}#1}}

\allowdisplaybreaks
\title[Uniqueness of continuous solutions of the Landau equation]{Weak--strong uniqueness for the Landau equation}
\author{Weinan Wang}
\address{Department of Mathematics, University of Oklahoma, Norman, OK 73019, USA}
\email{ww@ou.edu}
\date{\rev{November 19, 2026}}
\begin{document}

\begin{abstract}
In this paper, we consider the spatially inhomogeneous Landau equation with soft potentials,
including the Coulomb case. Let $f$ and $g$ be nonnegative, uniformly
continuous distributional solutions on $[0,T]$ with the same initial data.
We prove that $f=g$ on $[0,T]$ if both solutions have sufficient polynomial
decay in velocity and the diffusion matrix of $f$ satisfies the usual
anisotropic lower bound. For weak solutions, this proves a stronger form of
\cite[Conjecture~1.4]{HW}, with no prescribed modulus of continuity for the
initial datum. The proof is by duality. We obtain a weighted $L^2$ estimate
for the velocity Hessian of the solution to a backward adjoint equation.
In the second-order nonlocal term, we write $g$ as the sum of its velocity
mollification and a remainder. Uniform continuity makes the remainder small.
For the mollified part, two integrations by parts give a bound proportional
to the length of the time interval.
\end{abstract}\maketitle

\section{Introduction}

In this paper, we study uniqueness for the Cauchy problem for the spatially inhomogeneous Landau equation with soft potentials in three dimensions. Henderson, Snelson and Tarfulea \cite{HST} constructed solutions for bounded initial densities with polynomial decay in velocity and a quantitative lower bound. Their solutions are classical at positive times. For uniqueness, the difficulty is that the velocity Hessian may be unbounded as $t\searrow0$, even for uniformly continuous initial data.

Later, in \cite[Theorem 1.2]{HW}, Henderson and the author of this paper  proved weak--strong uniqueness for initial data with sufficiently fast velocity decay, $C^{\alpha/3}$ regularity in $x$, and a velocity modulus of order $|\log r|^{-\theta}$, provided $\theta\alpha/[2(2+\alpha)]>1$. They conjectured that for $\theta>2$ one could dispense with the spatial H\"older assumption \cite[Conjecture 1.4]{HW}. 

\subsection{The equation and the solution classes}\label{sec:equation}
Write $Y=\partial_t+v\cdot\nabla_x$, where $(t,x,v)\in[0,T]\times\R^3\times\R^3$. We take $\gamma\in[-3,0)$ and set, for $z\ne0$ and a fixed $a_\gamma>0$,
\[
a(z)=a_\gamma|z|^{\gamma+2}
 \left(\mathrm{Id}-\frac{z\otimes z}{|z|^2}\right),
\qquad b_i(z)=\partial_j a_{ij}(z)=-2a_\gamma |z|^\gamma z_i.
\]
Repeated velocity indices are summed. At each $(t,x)$ we define
\[
 \ba^h=a*_v h,\qquad \beta^h=b*_v h,\qquad
 \bc^h=
 \begin{cases}
  c_\gamma|\cdot|^\gamma *_v h,&-3<\gamma<0,\\
  c_{-3}h,&\gamma=-3.
 \end{cases}
\]
Here $c_\gamma=2a_\gamma(\gamma+3)$ for $\gamma>-3$ and $c_{-3}=8\pi a_{-3}$. We write $\jv=(1+|v|^2)^{1/2}$ and $\|h\|_{L^{\infty,k}}=\|\jv^kh\|_{L^\infty}$. When $h\in L^{\infty,\ell}$ with $\ell>5+\gamma$, the coefficients are locally bounded. In the sense of distributions, we have
\begin{equation}\label{div}
 \partial_i\ba^h_{ij}=\beta^h_j,\qquad \partial_j\beta^h_j=-\bc^h .
\end{equation}
For smooth functions these identities give
\begin{equation}\label{L}
 \begin{split}
 Yh&=Q(h,h),\\
 Q(h,u)&=\ba^h:D_v^2 u+\bc^h u
       =\nabla_v\cdot(\ba^h\nabla_vu-\beta^h u).
 \end{split}
\end{equation}

\Needspace{4\baselineskip}
\begin{definition*}
Let $h\in L^{\infty,\ell}((0,T)\times\R^6)$, with $\ell>5+\gamma$.
\begin{itemize}
\item[(i)] We call $h$ a \emph{distributional solution} of \eqref{L} if, for every $\varphi\in C_c^\infty((0,T)\times\R^6)$,
\begin{equation}\label{dist}
 \iiint h\big(Y\varphi+\ba^h:D_v^2\varphi+2\beta^h\cdot\nabla_v\varphi\big)\dd v\dd x\dd t=0.
\end{equation}
\item[(ii)] We call $h$ a \emph{weak solution} if $\nabla_vh\in L^2_{\rm loc}((0,T)\times\R^6)$ and, for the same test functions,
\begin{equation}\label{weakform}
 \iiint h\,Y\varphi\dd v\dd x\dd t
 =\iiint(\ba^h\nabla_vh-\beta^h h)\cdot\nabla_v\varphi\dd v\dd x\dd t .
\end{equation}
\end{itemize}
\end{definition*}

Every weak solution in (ii) satisfies (i) by integration by parts in $v$ (Lemma~\ref{lem:dual}). The kinetic weak solutions in \cite[Section 5]{HST} satisfy (ii), since the kinetic energy space gives a locally square-integrable velocity gradient and \rev{admits smooth compactly supported test functions}. Their drift $\bar b^h$ equals $-\beta^h$ in our notation. We also use the class $C^2_{\rm kin,loc}$ of functions for which $h$, $\nabla_vh$, $D_v^2h$, and $Yh$ are continuous. For a pointwise solution in this class, integration along $t\mapsto(t,y+tv,v)$ and integration by parts in $v$ give \eqref{weakform}.

\subsection{Main results}

\begin{theorem}\label{thm:main}
Fix $T,K,\lambda>0$ and $k\ge44$. Suppose that $f,g:[0,T]\times\R^6\to[0,\infty)$ satisfy the following hypotheses.
\begin{itemize}
\item[(A1)] Both functions are uniformly continuous on the closed time slab and solve \eqref{L} in the sense of distributions on its interior. They have the same initial values, and
\[
 \|f\|_{L^{\infty,k}}+\|g\|_{L^{\infty,k}}\le K.
\]
\item[(A2)] For every $(t,x,v)\in[0,T]\times\R^6$ and $e\in\mathbb S^2$,
\[
 e\cdot\ba^f(t,x,v)e\ge\lambda\jv^\gamma,\qquad
 e\cdot\ba^f(t,x,v)e\ge\lambda\jv^{\gamma+2}\quad\text{if }e\cdot v=0.
\]
\end{itemize}
Then $f=g$ on $[0,T]\times\R^6$.
\end{theorem}

The powers in (A2) differ by two because the radial diffusion is weaker at large $|v|$ than the diffusion in directions orthogonal to $v$. For the solutions constructed in \cite{HST}, (A2) follows from the positivity condition \cite[(1.7)]{HW} on the initial density \cite[Lemma 2.5 and the ensuing remark]{HST}. \rev{That condition requires the initial density to have, at every $x$, a common positive lower bound on a velocity ball of fixed radius whose center is uniformly bounded.}

\begin{corollary}[Weak--strong uniqueness]\label{cor}
Fix $k\ge44$. Let $f_{\rm in}\in L^{\infty,k}(\R^6)$ satisfy the nondegeneracy condition \cite[(1.7)]{HW}, and let $f\in L^{\infty,k}([0,T]\times\R^6)$ be a solution with initial datum $f_{\rm in}$ constructed in \cite[Theorem 1.2]{HST}. Suppose that $g:[0,T]\times\R^6\to[0,\infty)$ is uniformly continuous, is a weak solution of \eqref{L} on $(0,T)\times\R^6$, and belongs to $L^{\infty,5+\gamma+\eta}([0,T]\times\R^6)$ for some $\eta>0$. If $g(t,x,v)\to f_{\rm in}(x,v)$ as $t\searrow0$ for every $(x,v)$, then $f=g$ on $[0,T]\times\R^6$.
\end{corollary}

\begin{remark}[Comparison with \cite{HW}]\label{rem:HW}
\rev{Since $g$ is uniformly continuous up to $t=0$, its pointwise initial trace is $g(0)=f_{\rm in}$. The hypotheses of Corollary~\ref{cor} therefore force $f_{\rm in}$ to be uniformly continuous in $(x,v)$, but they prescribe no modulus. For sufficiently large finite $k$, \cite[Theorem 1.2]{HW} gives uniqueness on an initial subinterval $[0,T_1]$, whereas Corollary~\ref{cor} gives uniqueness on all of $[0,T]$ for every $k\ge44$.}
\end{remark}

\begin{remark}[Nonnegativity of $g$]\label{rem:sign}
\rev{The assumption $g\ge0$ enters the proof of Corollary~\ref{cor} through Lemma~\ref{lem:decay}, which propagates the velocity decay of $g(0)$ by a comparison argument and requires $\ba^g\ge0$. For solutions in the weak class of \cite{HST} that belong to $L^{\infty,5+\eta'}$ with $\eta'>0$, nonnegativity follows from the proof of \cite[Lemma 5.1]{HST}.}
\end{remark}

\subsection{\texorpdfstring{\rev{Strategy of the proof}}{Strategy of the proof}}\label{sec:strategy}
Write $w=f-g$. In the expansion $Yw=Q(g,w)+Q(w,f)$, the term $Q(w,f)$ contains $\ba^w:D_v^2f$. To close a Gronwall estimate for $w$, Henderson and the author of this paper proved a bound of the form \cite[Proposition 3.1]{HW}
\[
 \int_0^{T_0}\|\jv^7 D_v^2f(t)\|_{L^\infty(\R^6)}\dd t<\infty .
\]
Kinetic scaling gives a factor $t^{-1}$ for two velocity derivatives of a bounded solution. The modulus of the initial data in \cite{HW} supplies the integrable logarithmic correction. Conjecture~1.3 of that paper proposed a Schauder estimate involving only velocity H\"older seminorms, which would have allowed the same comparison for data with less spatial regularity.

In \cite{DW}, Dong and the author of this paper recently disproved this Schauder estimate even for the constant-coefficient Kolmogorov equation. Their forcing $\cos(Nx_1)$ is constant in $v$, yet one entry of the velocity Hessian of the solution varies by order one across a velocity interval of length comparable to $N^{-1/3}$. Its $C^\alpha_v$ seminorm grows at least as $N^{\alpha/3}$.

\rev{We replace the pointwise bound for $D_v^2f$ by an $L^2$ bound for the velocity Hessian of a solution of an adjoint equation.} To derive the adjoint equation, write $Yw=Q(f,w)+Q(w,g)$ and integrate against a test function $\phi$. Formally,
\begin{equation}\label{formal}
 \frac{\dd}{\dd t}\iint w\phi\dd x\dd v
 =\iint w(\cP\phi+\K_g\phi)\dd x\dd v ,
\end{equation}
where
\begin{align*}
 \cP\phi&=Y\phi+\ba^f:D_v^2\phi+2\beta^f\cdot\nabla_v\phi,\\
 \K_g\phi&=a_{ij}*_v(g\,\partial_{ij}\phi)-2b_i*_v(g\,\partial_i\phi).
\end{align*}
The oddness of $b$ accounts for the minus sign in the last convolution. Neither operator differentiates a density. Lemma~\ref{lem:weak} justifies \eqref{formal} for solutions of the backward equation.

\rev{On an interval $I=(t_1,t_2)$ of length $h$ we solve the backward problem $\cP\phi=\Psi$, $\phi(t_2)=0$, for smooth compactly supported $\Psi$. The kinetic $L^2$ estimate of Dong and Yastrzhembskiy \cite{DY} requires uniform ellipticity, whereas the ellipticity ratio of $\ba^f$ grows with $|v|$. We therefore work in the rescaled boxes of \cite{CSS,HST}. There the diffusion matrices are uniformly elliptic and share a modulus of continuity in $(x,v)$, so the estimate of \cite{DY} bounds the Hessian with a constant independent of the central velocity. Summing over the boxes produces lower-order terms, which a weighted energy estimate controls. Since the densities may have no decay in $x$, the test functions must carry a weight that grows in $x$; we take it constant along the characteristics of $Y$, so that it commutes with the transport operator.}

\rev{The bounds for $\phi$ and $\nabla_v\phi$ gain factors $h$ and $h^{1/2}$, respectively, but the Hessian bound has no small factor as $h\searrow0$. The smallness of the part of $\K_g\phi$ that contains $D_v^2\phi$ must therefore come from $g$. Let $g_\epsilon$ be a mollification of $g$ in velocity. By uniform continuity and polynomial decay, $g-g_\epsilon$ is small in a weighted $L^\infty$ norm, while in the convolution containing $g_\epsilon D_v^2\phi$ two integrations by parts move both derivatives off $\phi$. In the weighted space $\cZ(I)$ of Section~\ref{sec:proof}, this gives
\[
 \|\K_g\phi\|_{\cZ(I)}
 \le \big(C\varpi_g(\epsilon)^{1/2}+C_\epsilon h+C h^{1/2}\big)
      \|\Psi\|_{\cZ(I)},
\]
where $\varpi_g$ is a uniform modulus of $g$ in $v$. We choose $\epsilon$ so that the first term in parentheses is at most $1/4$, and then $h$ so small that the other two terms sum to at most $1/4$. If $w(t_1)=0$, integrating \eqref{formal} then gives
\[
 \left|\iiint_{I\times\R^6} w\Psi\right|
 \le \tfrac12\,\|w\|_{\cZ(I)'}\|\Psi\|_{\cZ(I)}.
\]
Taking the supremum over smooth $\Psi$ with $\|\Psi\|_{\cZ(I)}\le1$ shows that $\|w\|_{\cZ(I)'}\le\frac12\|w\|_{\cZ(I)'}$, and hence $w=0$ on $I$. Since the choice of $h$ does not depend on $t_1$, iterating over consecutive intervals gives uniqueness up to $T$.}

\subsection{Related work}
\rev{The Landau equation arises from the Boltzmann equation in the grazing-collision limit; see Villani \cite{VillaniH} for the homogeneous setting and Alexandre and Villani \cite{AV} for the inhomogeneous one. The theory of the inhomogeneous equation is most complete for small perturbations. In a periodic box, existence and convergence to equilibrium near a Maxwellian were proved in \cite{Guo,SG}, weighted semigroup estimates are available for moderately soft potentials \cite{CTW} and for very soft and Coulomb potentials \cite{CM}, and global mild solutions with summable spatial Fourier coefficients were constructed in \cite{DLSS}. In the whole space, Luk \cite{Luk} uses dispersion to prove global stability near vacuum for localized data when $-2<\gamma<0$. We impose no smallness condition of either kind.}

\rev{As for large data, the construction in \cite{HST} relies on the H\"older estimate of Golse, Imbert, Mouhot and Vasseur \cite{GIMV}, on the Schauder estimates of Henderson and Snelson \cite{HS}, on the change of variables of Cameron, Silvestre and Snelson \cite{CSS}, and on the lower bounds obtained from the spreading of positivity in \cite{HSTmass}. We use the last two ingredients as well: the change of variables defines the boxes of Section~\ref{sec:frames}, and the lower bounds are the source of hypothesis {\rm(A2)}. The regularity estimates, in contrast, are interior estimates and do not control the velocity Hessian up to $t=0$.}

\rev{For this reason uniqueness has required more regularity than existence. The uniqueness theorem of \cite{HST} assumes H\"older regularity of the datum in both $x$ and $v$, and \cite{HW} relaxes this to H\"older regularity in $x$ and a logarithmic modulus in $v$. Both proofs rest on a bound for $D_v^2f$ that is integrable in time near $t=0$ and is derived from Schauder estimates. By \cite{DW}, the velocity-only Schauder estimate conjectured in \cite{HW}, which would have removed the spatial H\"older assumption, is false.}

\rev{A different approach is due to Alonso, Gualdani and Sun \cite{AGS}, who prove uniqueness of continuous solutions on $\mathbb T^3\times\R^3$ for a kinetic model and for the Landau--Coulomb equation with an additional viscosity $\nu\Delta_v$, $\nu>0$. They compare the densities after applying an operator of negative order in $x$ and $v$ and, as we do, obtain smallness from continuity through a mollification. Theorem~\ref{thm:main} treats \eqref{L} without viscosity. The lower bound {\rm(A2)} on one solution takes over the role of the added diffusion, the spatial domain is $\R^3$, and all $\gamma\in[-3,0)$ are allowed. The cost appears in the velocity weights: the estimate for the adjoint Hessian in Proposition~\ref{prop:solv} loses $(\gamma+2)_+$ powers of $\jv$. In addition, Kim, Guo and Hwang \cite{KGH} prove uniqueness of small $L^\infty$ perturbations of a Maxwellian in the Coulomb case by a weighted $L^2$ argument and kinetic $S_p$ estimates. We also combine weighted $L^2$ bounds with a kinetic Sobolev estimate, but we apply them to the adjoint equation and assume no smallness.}

\rev{The spatially homogeneous equation is better understood. Uniqueness has been proved there by Wasserstein and coupling arguments \cite{FG,FH}, for densities with an integrable-in-time $L^\infty_v$ norm at the Coulomb endpoint \cite{F}, in $L^\infty_tL^p_v$ with $p>3/2$ \cite{CG}, by relative entropy \cite{Tabary}, and in a weighted $H^{-1/2}$ space \cite{HJL}. As for regularity, the Fisher information is nonincreasing, which propagates smoothness globally \cite{GS}. In the Coulomb case, entropy dissipation controls a weighted $L^2$ norm of $\nabla_v\sqrt f$ \cite{DesvD}, and global smooth solutions exist for weighted $L^{3/2}$ data under additional moment and lower bound assumptions \cite{GGL}. None of these arguments has to deal with the transport operator $Y$.}

\rev{Turning to the linear theory, the gain of regularity in $x$ comes from the interplay of velocity diffusion and transport. Kolmogorov's fundamental solution \cite{Kolmogorov} exhibits the scales $r$, $r^2$ and $r^3$ in velocity, time and space, in a frame moving with the central velocity; these are the scales of the map $\Lambda$ in the proof of Proposition~\ref{prop:local}. The commutator $[\partial_{v_i},Y]=\partial_{x_i}$ gives the bracket condition of H\"ormander \cite{Hormander}, and Bouchut \cite{Bouchut} quantified the gain in spatial Sobolev regularity. Local $L^p$ estimates for ultraparabolic operators with coefficients of vanishing mean oscillation were obtained by Bramanti, Cerutti and Manfredini \cite{BCM}. For the backward problem we use the global nondivergence-form theory of Dong and Yastrzhembskiy \cite{DY}, which allows coefficients that are merely measurable in time and gives both the Hessian estimate and solvability with zero initial data. We need the nondivergence form because $\K_g\phi$ contains the full velocity Hessian of $\phi$; for divergence-form equations we refer to the variational theory of Albritton, Armstrong, Mourrat and Novack \cite{AAMN} and to the estimates of \cite{DYDiv}. In addition, for interior Schauder estimates we refer to Manfredini \cite{Manfredini} and to Imbert and Mouhot \cite{IM}, and for coefficients measurable in time to \cite{BiagiBramanti,HW,DYSchauder}; see also \cite{PLP}.}

\section{Frames, boxes, weights}\label{sec:frames}

\subsection{Anisotropic frames and boxes}
We use the anisotropic rescaling of \cite{CSS,HST}. Set
\[
\kappa:=\Big(1+\frac{\gamma}{2}\Big)_+\in[0,1),\qquad\text{so that}\qquad 2\kappa=(\gamma+2)_+ .
\]
For $v\ne0$ let $P_v$ be the orthogonal projection onto $\R v$. Set $P_0:=0$ and $P_v^\perp:=\mathrm{Id}-P_v$. We define
\[
S_v:=\jv^{\gamma/2}P_v+\jv^{1+\gamma/2}P_v^\perp\quad(S_0:=\mathrm{Id}),\qquad r_v:=\jv^{-\kappa},\qquad
E_\rho(y):=\{y+r_yS_y\tv:\ |\tv|<\rho\}.
\]
The matrix $S_v$ is symmetric and positive definite, and for every $\xi\in\R^3$
\begin{equation}\label{Sbasic}
|S_v\xi|^2=\jv^{\gamma}|P_v\xi|^2+\jv^{\gamma+2}|P_v^\perp\xi|^2,\qquad
\jv^{\gamma/2}|\xi|\le|S_v\xi|\le\jv^{1+\gamma/2}|\xi|,\qquad
r_v|S_v\xi|\le|\xi| .
\end{equation}
The last inequality holds because $1+\gamma/2\le\kappa$. It implies that $E_\rho(y)\subset B_\rho(y)$.

 The letters $C$ and $c$ denote positive constants which may change from line to line. They
depend only on $\gamma$, $T$, $K$, $k$, $\lambda$ and on the modulus of continuity of $f$ (and of $g$ where indicated).

\begin{lemma}[Frames]\label{lem:frames}
There is a constant $C_0\ge1$ such that for all $v,y\in\R^3$ with $|v-y|\le4$ and all $\xi\in\R^3$,
\begin{equation}\label{comp}
C_0^{-1}|S_y\xi|\le|S_v\xi|\le C_0|S_y\xi|,\qquad C_0^{-1}r_y\le r_v\le C_0r_y .
\end{equation}
Moreover, there are constants $0<c<C$ such that for every measurable function $G\ge0$ on $\R^3$,
\begin{equation}\label{aver}
c\int_{\R^3}G\dd v\le\int_{\R^3}\Big(\int_{E_{1/2}(y)}G\dd v\Big)\frac{\dd y}{\det(r_yS_y)},\qquad
\int_{\R^3}\Big(\int_{E_{1}(y)}G\dd v\Big)\frac{\dd y}{\det(r_yS_y)}\le C\int_{\R^3}G\dd v .
\end{equation}
\end{lemma}

\begin{proof}
Let $|v-y|\le4$.

\step{Step 1: comparison of $r_v$ and $r_y$.}
Since $\jv\le\jy+|v-y|\le5\jy$ and, by symmetry, $\jy\le5\jv$, the second part of \eqref{comp} holds with the constant $5^\kappa$.

\step{Step 2: comparison of $S_v$ and $S_y$ for $|y|\le8$.}
In this case $|v|\le12$, and by \eqref{Sbasic} all eigenvalues of $S_y$ and $S_v$ lie between two positive constants which depend only
on $\gamma$. This gives the first part of \eqref{comp}.

\step{Step 3: comparison of $S_v$ and $S_y$ for $|y|>8$.}
Now $|v|\ge|y|-4\ge|y|/2$. Writing $\hat v:=v/|v|$ and $\hat y:=y/|y|$, we have
\[
|\hat v-\hat y|\le\frac{2|v-y|}{|y|}\le\frac{8}{|y|},\qquad\text{hence}\qquad
\|P_v-P_y\|\le2|\hat v-\hat y|\le\frac{16}{|y|} .
\]
Consequently
\[
|P_v^\perp\xi|\le|P_y^\perp\xi|+\frac{16}{|y|}\,|\xi| .
\]
We insert this bound in the first identity of \eqref{Sbasic} and use Step 1 and $|P_v\xi|\le|\xi|$:
\begin{align*}
|S_v\xi|^2&=\jv^{\gamma}|P_v\xi|^2+\jv^{\gamma+2}|P_v^\perp\xi|^2\\
&\le C\jy^{\gamma}|\xi|^2+C\jy^{\gamma+2}\Big(|P_y^\perp\xi|^2+\frac{|\xi|^2}{|y|^2}\Big)\\
&\le C\Big(\jy^{\gamma}|\xi|^2+\jy^{\gamma+2}|P_y^\perp\xi|^2\Big)\le C|S_y\xi|^2 .
\end{align*}
In the last step we used that both $\jy^{\gamma}|\xi|^2$ and $\jy^{\gamma+2}|P_y^\perp\xi|^2$ are bounded by $|S_y\xi|^2$, by
\eqref{Sbasic}. This is the upper bound in \eqref{comp}. The lower bound follows by exchanging the roles of $v$ and $y$: the argument
only used that $|y|>4$ and $|\hat v-\hat y|\le8/|y|$, and we also have $|v|>4$ and $|\hat v-\hat y|\le8/|v|$.

\step{Step 4: reduction of \eqref{aver} to bounds for a density.}
For $\rho\in\{\frac12,1\}$ define
\[
\mu_\rho(v):=\int_{\R^3}\mathbf 1_{E_\rho(y)}(v)\,\frac{\dd y}{\det(r_yS_y)} .
\]
By Fubini's theorem,
\[
\int_{\R^3}\Big(\int_{E_{\rho}(y)}G\dd v\Big)\frac{\dd y}{\det(r_yS_y)}=\int_{\R^3}G(v)\,\mu_\rho(v)\dd v .
\]
Therefore \eqref{aver} is equivalent to the two bounds $\mu_{1/2}\ge c$ and $\mu_1\le C$ on $\R^3$. We recall that
$|E_\rho(v)|=\det(r_vS_v)\,|B_\rho|$, and we note that \eqref{comp} implies, for $|v-y|\le4$ and $\eta\in\R^3$,
\begin{equation}\label{compinv}
\big|(r_vS_v)^{-1}\eta\big|\le C_0^2\,\big|(r_yS_y)^{-1}\eta\big|,\qquad C_0^{-6}\det(r_vS_v)\le\det(r_yS_y)\le C_0^{6}\det(r_vS_v).
\end{equation}
Indeed, the first inequality follows from \eqref{comp} by duality, since $|S^{-1}\eta|=\sup_{\xi\ne0}\eta\cdot\xi/|S\xi|$ for a symmetric
positive definite matrix $S$; the second one follows from the comparison of the quadratic forms $|S_v\xi|^2$ and $|S_y\xi|^2$ and of
$r_v$ and $r_y$.

\step{Step 5: the upper bound $\mu_1\le C$.}
Let $v\in E_1(y)$. Then $|v-y|<1$ and $|(r_yS_y)^{-1}(y-v)|<1$, so by \eqref{compinv}
\[
\big|(r_vS_v)^{-1}(y-v)\big|<C_0^2,\qquad\text{that is,}\qquad y\in E_{C_0^2}(v) .
\]
Thus the set $\{y:\ v\in E_1(y)\}$ is contained in $E_{C_0^2}(v)\cap B_1(v)$. Using the lower bound for $\det(r_yS_y)$ in
\eqref{compinv} on this set, we obtain
\[
\mu_1(v)=\int_{\{y:\,v\in E_1(y)\}}\frac{\dd y}{\det(r_yS_y)}\le C_0^6\,\frac{|E_{C_0^2}(v)|}{\det(r_vS_v)}=C_0^6\,|B_{C_0^2}|=:C .
\]

\step{Step 6: the lower bound $\mu_{1/2}\ge c$.}
Let $\rho_0:=1/(2C_0^2)$ and let $y\in E_{\rho_0}(v)$. Then $|y-v|<1$, and by \eqref{compinv} with the roles of $v$ and $y$ exchanged,
\[
\big|(r_yS_y)^{-1}(v-y)\big|\le C_0^2\,\big|(r_vS_v)^{-1}(v-y)\big|<C_0^2\rho_0=\frac12,\qquad\text{that is,}\qquad v\in E_{1/2}(y) .
\]
Using the upper bound for $\det(r_yS_y)$ in \eqref{compinv}, we conclude
\[
\mu_{1/2}(v)\ge\int_{E_{\rho_0}(v)}\frac{\dd y}{\det(r_yS_y)}\ge C_0^{-6}\,\frac{|E_{\rho_0}(v)|}{\det(r_vS_v)}=C_0^{-6}\,|B_{\rho_0}|=:c .
\qedhere
\]
\end{proof}

\subsection{Weights}
Let $n_x:=4$. For $t\in[0,T]$ and $\epsilon\in[0,1]$ we define
\begin{gather*}
\theta(t,x,v):=\big\langle (x-tv)\jv^{-\kappa}\big\rangle^{n_x},\qquad \te:=\frac{\theta}{1+\epsilon\theta}\quad(\theta_0=\theta),\\
\omega_{m',\epsilon}:=\te\,\jv^{-m'},\qquad\omega_{m'}:=\omega_{m',0}.
\end{gather*}
For $\epsilon>0$ the function $\te$ is bounded, and $\te$ increases to $\theta$ as $\epsilon\searrow0$. We fix once and for all the exponents
\[
m_1:=8,\qquad n:=8,\qquad m:=m_1+2\kappa,\qquad\text{so that}\qquad \omega_m=r_v^{2}\,\omega_{m_1}.
\]

\begin{lemma}[Weights]\label{lem:weights}
Let $t\in[0,T]$ and $\epsilon\in[0,1]$. The following holds with constants independent of $\epsilon$.
\begin{itemize}
\item[\rm(W1)] \emph{Transport.} $Y\te=0$.
\item[\rm(W2)] \emph{Derivatives in $v$.} We have
\[
|\nabla_v\te|\le C\jv^{-\kappa}\te,\qquad |D_v^2\te|\le C\jv^{-2\kappa}\te ,
\]
and consequently
\[
|S_v\nabla_v\te|\le C\te,\qquad |S_vD_v^2\te\, S_v|\le C\te .
\]
\rev{Moreover, for every $m'\ge0$, with a constant $C$ depending also on $m'$,}
\[
|S_v\nabla_v\omega_{m',\epsilon}|\le C\omega_{m',\epsilon},\qquad
|\nabla_v\omega_{m',\epsilon}|\le C\big(\jv^{-\kappa}+\jv^{-1}\big)\,\omega_{m',\epsilon} .
\]
\item[\rm(W3)] \emph{Comparison.} For all $x,v,u\in\R^3$,
\[
\te(t,x,v)\le C\langle v-u\rangle^{n_x}\,\te(t,x,u),\qquad \langle x\rangle^{n_x}\le C\jv^{2n_x}\,\theta(t,x,v) .
\]
\end{itemize}
\end{lemma}
\begin{proof}
We first prove {\rm(W1)--(W3)} for $\epsilon=0$, that is, for $\theta$, and then pass to $\te$. We write
\[
\theta=F(s),\qquad F(s):=\langle s\rangle^{n_x},\qquad s=s(t,x,v):=(x-tv)\,\jv^{-\kappa}\in\R^3,
\]
and note that
\begin{equation}\label{Fder}
|\nabla F(s)|\le C\langle s\rangle^{n_x-1},\qquad |D^2F(s)|\le C\langle s\rangle^{n_x-2}.
\end{equation}

\step{Step 1: {\rm(W1)} for $\theta$.}
Since $Y(x-tv)=-v+v=0$ and $\jv$ does not depend on $(t,x)$, we have $Ys=0$ and therefore $Y\theta=0$.

\step{Step 2: {\rm(W2)} for $\theta$.}
Differentiating $s$ with respect to $v_i$ gives
\[
\partial_{v_i}s=-t\,\jv^{-\kappa}e_i-\kappa\,s\,v_i\,\jv^{-2} .
\]
Using $t\le T$ and $\jv^{-1}\le\jv^{-\kappa}$, we obtain
\[
|\partial_{v_i}s|\le C\langle s\rangle\jv^{-\kappa},\qquad
|\partial_{v_iv_j}s|\le C\langle s\rangle\jv^{-\kappa-1}\le C\langle s\rangle\jv^{-2\kappa}.
\]
By the chain rule and \eqref{Fder},
\begin{align*}
|\nabla_v\theta|&\le|\nabla F(s)|\,|\nabla_vs|\le C\langle s\rangle^{n_x}\jv^{-\kappa}=C\jv^{-\kappa}\theta,\\
|D_v^2\theta|&\le|D^2F(s)|\,|\nabla_vs|^2+|\nabla F(s)|\,|D_v^2s|\le C\langle s\rangle^{n_x}\jv^{-2\kappa}=C\jv^{-2\kappa}\theta .
\end{align*}
By \eqref{Sbasic} the operator norm of $S_v$ is $\jv^{1+\gamma/2}\le\jv^{\kappa}$, so these two bounds imply
\[
|S_v\nabla_v\theta|\le C\theta,\qquad |S_vD_v^2\theta\,S_v|\le C\theta .
\]
It remains to treat the factor $\jv^{-m'}$. Its gradient $\nabla_v\jv^{-m'}=-m'\jv^{-m'-2}v$ is parallel to $v$, so that
\[
\big|S_v\nabla_v\jv^{-m'}\big|=m'\jv^{-m'-2}\jv^{\gamma/2}|v|\le m'\jv^{-m'},\qquad
\big|\nabla_v\jv^{-m'}\big|\le m'\jv^{-m'-1}.
\]
The product rule applied to $\omega_{m'}=\theta\jv^{-m'}$ now gives the last two bounds of {\rm(W2)}.

\step{Step 3: {\rm(W3)} for $\theta$.}
Let $s_v:=(x-tv)\jv^{-\kappa}$ and $s_u:=(x-tu)\langle u\rangle^{-\kappa}$. We use the two inequalities
\[
|x-tv|\le|x-tu|+T|v-u|,\qquad \jv^{-\kappa}\le C\langle u\rangle^{-\kappa}\langle v-u\rangle^{\kappa}.
\]
They give
\[
|s_v|\le C\,|x-tu|\,\langle u\rangle^{-\kappa}\langle v-u\rangle^{\kappa}+T|v-u|\le C\langle s_u\rangle\langle v-u\rangle,
\]
hence $\langle s_v\rangle\le C\langle s_u\rangle\langle v-u\rangle$. Raising this to the power $n_x$ yields the first inequality of
{\rm(W3)}. For the second one we observe that
\[
\langle x\rangle\le C\langle x-tv\rangle\jv\qquad\text{and}\qquad\langle x-tv\rangle\le\langle s_v\rangle\jv^{\kappa},
\]
so that $\langle x\rangle\le C\langle s_v\rangle\jv^{1+\kappa}\le C\langle s_v\rangle\jv^{2}$, and we raise to the power $n_x$ again.

\step{Step 4: from $\theta$ to $\te$.}
We have $\te=F_\epsilon(\theta)$ with $F_\epsilon(\sigma):=\sigma/(1+\epsilon\sigma)$ for $\sigma\ge1$. This function is increasing and
satisfies
\[
0\le F_\epsilon'(\sigma)\le\frac{F_\epsilon(\sigma)}{\sigma},\qquad
|F_\epsilon''(\sigma)|\le\frac{2F_\epsilon(\sigma)}{\sigma^2},\qquad
F_\epsilon(C\sigma)\le CF_\epsilon(\sigma)\quad(C\ge1).
\]
Consequently $Y\te=F_\epsilon'(\theta)\,Y\theta=0$, and
\[
|\nabla_v\te|\le\frac{\te}{\theta}\,|\nabla_v\theta|,\qquad
|D_v^2\te|\le\frac{\te}{\theta}\,|D_v^2\theta|+\frac{2\te}{\theta^2}\,|\nabla_v\theta|^2 .
\]
Thus the bounds of {\rm(W2)} for $\theta$ imply the same bounds for $\te$, with constants independent of $\epsilon$, and the bounds
for $\omega_{m',\epsilon}$ follow as in Step 2. Finally, the first inequality of {\rm(W3)} for $\te$ follows from the one for $\theta$,
because $F_\epsilon$ is increasing and $F_\epsilon(C\sigma)\le CF_\epsilon(\sigma)$.
\end{proof}

\subsection{Velocity convolutions}
\begin{lemma}\label{lem:int}
Let $q\ge9$. For all $v\in\R^3$ we have
\begin{align}
&\int_{\R^3}|v-u|^{s}\langle u\rangle^{-q}\dd u\le C\jv^{s}&&\text{for }-3<s\le2,\label{I1}\\
&\int_{\R^3}\min\big\{|v-u|^{s+j},\,|v-u|^{s}\langle u\rangle^{j}\big\}\langle u\rangle^{-q}\dd u\le C\jv^{s}&&\text{for }-3\le s<0,\ j\in\{1,2\}.\label{I2}
\end{align}
Here the constants may depend on $q$ and $s$.
\end{lemma}
\begin{proof}
We split the domain of integration into the regions $\{|v-u|<1\}$ and $\{|v-u|\ge1\}$.

\step{Step 1: the region $|v-u|<1$.}
Here $\langle u\rangle^{-q}\le C\jv^{-q}$. In \eqref{I1} the remaining factor $|v-u|^{s}$ is integrable over the unit ball because
$s>-3$. In \eqref{I2} we bound the minimum by its first entry $|v-u|^{s+j}$, which is integrable because $s+j\ge-2$. In both cases
this region contributes at most $C\jv^{-q}\le C\jv^{s}$.

\step{Step 2: the region $|v-u|\ge1$ for $s<0$.}
Here $|v-u|^{s}\le C\langle v-u\rangle^{s}$, and in \eqref{I2} we bound the minimum by its second entry. Thus it suffices to show
that, for $j\in\{0,1,2\}$,
\begin{equation}\label{I3}
\int_{\R^3}\langle v-u\rangle^{s}\langle u\rangle^{j-q}\dd u\le C\jv^{s}\qquad(-3\le s<0).
\end{equation}
On the set $\{|u|\le|v|/2\}$ we have $\langle v-u\rangle\ge\jv/2$, hence $\langle v-u\rangle^{s}\le C\jv^{s}$, and
$\langle u\rangle^{j-q}$ is integrable. On the set $\{|u|>|v|/2\}$ we have $\langle u\rangle\ge\jv/2$, hence
\[
\langle u\rangle^{j-q}\le C\jv^{-3}\langle u\rangle^{j-q+3},
\]
which is integrable because $q-j-3>3$; since $\langle v-u\rangle^{s}\le1$, this set contributes at most $C\jv^{-3}\le C\jv^{s}$.
This proves \eqref{I3}.

\step{Step 3: the region $|v-u|\ge1$ for $0\le s\le2$.}
This case only occurs in \eqref{I1}. We use $|v-u|^{s}\le C(\jv^{s}+\langle u\rangle^{s})$ and the integrability of
$\langle u\rangle^{s-q}$.
\end{proof}

\section{The coefficients}\label{sec:coef}
In this section $f$ satisfies
\[
0\le f\le K\jv^{-k}\quad\text{on }[0,T]\times\R^6,\qquad k\ge18,
\]
together with {\rm(A2)}. We abbreviate $\ba=\ba^f$, $\beta=\beta^f$, $\bc=\bc^f$, and we denote by $\varpi$ a modulus of continuity of
$f$ in $x$:
\[
|f(t,x,v)-f(t,x',v)|\le\varpi(|x-x'|),\qquad\varpi\ \text{nondecreasing},\quad\varpi(0+)=0 .
\]

For $z\ne0$ and $\xi,\eta\in\R^3$,
\begin{equation}\label{wedge}
\xi\cdot a(z)\eta=a_\gamma|z|^{\gamma}\,(z\wedge\xi)\cdot(z\wedge\eta).
\end{equation}
Indeed, by the definition of $a$ and Lagrange's identity,
\[
\xi\cdot a(z)\eta=a_\gamma|z|^{\gamma}\big(|z|^2\,\xi\cdot\eta-(z\cdot\xi)(z\cdot\eta)\big)=a_\gamma|z|^{\gamma}\,(z\wedge\xi)\cdot(z\wedge\eta).
\]
In particular $a(z)$ is symmetric and nonnegative.

\begin{lemma}\label{lem:coef}
The coefficients $\ba$, $\beta$, $\bc$ have the following properties.
\begin{itemize}
\item[\rm(a)] \emph{Size and regularity.} We have
\[
|\ba|\le C\jv^{2\kappa},\qquad |\beta|\le C\jv^{(\gamma+1)_+},\qquad 0\le\bc\le C ,
\]
and $\ba$ is continuously differentiable in $v$, with $\partial_i\ba_{ij}=\beta_j$.
\item[\rm(b)] \emph{Anisotropic ellipticity.} There are constants $0<c_1\le C_1$ such that for all $(t,x,v)$ and $\xi\in\R^3$,
\[
c_1|S_v\xi|^2\le\xi\cdot\ba(t,x,v)\xi\le C_1|S_v\xi|^2 .
\]
\item[\rm(c)] \emph{Continuity in normalized form.} For every $y\in\R^3$, all $v,v'\in E_2(y)$, all $x,x'\in\R^3$ and all $t\in[0,T]$,
\[
\big|S_y^{-1}\big(\ba(t,x,v)-\ba(t,x',v')\big)S_y^{-1}\big|\le C\Big(\big|(r_yS_y)^{-1}(v-v')\big|+\varpi(|x-x'|)^{1/2}\Big).
\]
\item[\rm(d)] \emph{Normalized drift.} $|r_vS_v^{-1}\beta(t,x,v)|\le C$.
\end{itemize}
\end{lemma}

\begin{proof}[Proof of {\rm(a)}]
By \eqref{I1} with $s=\gamma+2$ and with $s=\gamma+1$,
\begin{align*}
|\ba(t,x,v)|&\le CK\int_{\R^3}|v-u|^{\gamma+2}\langle u\rangle^{-k}\dd u\le C\jv^{\gamma+2}\le C\jv^{2\kappa},\\
|\beta(t,x,v)|&\le CK\int_{\R^3}|v-u|^{\gamma+1}\langle u\rangle^{-k}\dd u\le C\jv^{\gamma+1}\le C\jv^{(\gamma+1)_+}.
\end{align*}
If $\gamma>-3$, the bound for $\bc$ is \eqref{I1} with $s=\gamma$. If $\gamma=-3$, it holds because $\bc=c_{-3}f$. Finally,
$|\nabla a(z)|\le C|z|^{\gamma+1}$ is locally integrable and $f$ is bounded and decays. Hence $\ba$ is continuously differentiable in
$v$ with $\nabla_v\ba=(\nabla a)*_vf$, and in particular $\partial_i\ba_{ij}=b_j*_vf=\beta_j$.
\end{proof}

\begin{proof}[Proof of {\rm(b)}]
We fix $(t,x)$ and omit these variables.

\step{Step 1: the upper bound.}
Let $|v|\ge1$. We decompose
\[
\xi=\xi_r\hat v+\xi_\perp,\qquad\xi_\perp\perp v,
\]
and we write $z=v-u$. Since $z\wedge\hat v=-u\wedge\hat v$, we have $|z\wedge\hat v|=|P_v^\perp z|=|P_v^\perp u|\le\min\{|z|,|u|\}$, and
therefore
\[
|z\wedge\xi|\le|\xi_r|\,|z\wedge\hat v|+|z|\,|\xi_\perp|\le|\xi_r|\min\{|z|,|u|\}+|z|\,|\xi_\perp| .
\]
By \eqref{wedge},
\[
\xi\cdot a(z)\xi=a_\gamma|z|^{\gamma}|z\wedge\xi|^2\le2a_\gamma\Big(\xi_r^2\min\big\{|z|^{\gamma+2},|z|^{\gamma}|u|^2\big\}+|\xi_\perp|^2|z|^{\gamma+2}\Big).
\]
We multiply by $f(u)\le K\langle u\rangle^{-k}$ and integrate in $u$. By \eqref{I2} with $(s,j)=(\gamma,2)$ and \eqref{I1} with
$s=\gamma+2$,
\[
\xi\cdot\ba(v)\xi\le CK\big(\jv^{\gamma}\xi_r^2+\jv^{\gamma+2}|\xi_\perp|^2\big)=CK|S_v\xi|^2 .
\]
For $|v|<1$ the upper bound follows from {\rm(a)}, since then $|S_v\xi|\ge c|\xi|$.

The same computation shows that, for every measurable function $H$ on
$\R^3$ with $0\le H\le K'\langle u\rangle^{-9}$,
\begin{equation}\label{upperH}
\xi\cdot\ba^{H}(v)\xi\le CK'|S_v\xi|^2\qquad\text{for all }v,\xi\in\R^3 .
\end{equation}

\step{Step 2: the lower bound.}
For $|v|<1$ the lower bound follows from {\rm(A2)}, because $|S_v\xi|\le C|\xi|$. Let $|v|\ge1$. Since $\ba$ is symmetric and
nonnegative,
\[
\|\xi\|_a:=(\xi\cdot\ba(v)\xi)^{1/2}
\]
is a seminorm on $\R^3$. Let $\zeta:=S_v\xi$, and decompose $\zeta=\zeta_r\hat v+\zeta_\perp$ with $\zeta_\perp\perp v$. Then
\[
\xi=S_v^{-1}\zeta=\jv^{-\gamma/2}\zeta_r\hat v+\jv^{-1-\gamma/2}\zeta_\perp .
\]
By {\rm(A2)} for arbitrary directions, and since $|\xi|^2\ge\jv^{-\gamma}\zeta_r^2$,
\[
\|\xi\|_a^2\ge\lambda\jv^{\gamma}|\xi|^2\ge\lambda\zeta_r^2 .
\]
By {\rm(A2)} for directions orthogonal to $v$, the triangle inequality for $\|\cdot\|_a$, and the upper bound of Step~1,
\[
\sqrt\lambda\,|\zeta_\perp|\le\big\|\jv^{-1-\gamma/2}\zeta_\perp\big\|_a\le\|\xi\|_a+\big\|\jv^{-\gamma/2}\zeta_r\hat v\big\|_a
\le\|\xi\|_a+\sqrt{C_1}\,|\zeta_r| .
\]
Combining the two bounds, we find
\[
|S_v\xi|=|\zeta|\le|\zeta_r|+|\zeta_\perp|\le\lambda^{-1/2}\Big(2+\sqrt{C_1/\lambda}\Big)\|\xi\|_a .
\]
This is the lower bound, with $c_1:=\lambda\big(2+\sqrt{C_1/\lambda}\big)^{-2}$.
\end{proof}

\begin{proof}[Proof of {\rm(c)}]
We treat the dependence on $x$ and on $v$ separately; {\rm(c)} then follows from the triangle inequality, applied to
\[
\ba(t,x,v)-\ba(t,x',v')=\big[\ba(t,x,v)-\ba(t,x',v)\big]+\big[\ba(t,x',v)-\ba(t,x',v')\big].
\]

\step{Step 1: dependence on $x$.}
Fix $t$, $x$, $x'$ and put
\[
h:=f(t,x,\cdot)-f(t,x',\cdot),\qquad\mu:=\varpi(|x-x'|).
\]
Then $|h|\le\mu$ and $|h|\le2K\langle u\rangle^{-k}$. Taking the geometric mean of these two bounds, we find
\[
|h(u)|\le H(u):=(2K\mu)^{1/2}\langle u\rangle^{-k/2},\qquad\frac k2\ge9 .
\]
By \eqref{wedge} and the Cauchy--Schwarz inequality,
\[
|\xi\cdot\ba^{h}(v)\eta|\le\int_{\R^3}a_\gamma|z|^{\gamma}\,|z\wedge\xi|\,|z\wedge\eta|\,H(u)\dd u
\le\big(\xi\cdot\ba^{H}(v)\xi\big)^{1/2}\big(\eta\cdot\ba^{H}(v)\eta\big)^{1/2},
\]
where $z=v-u$. By \eqref{upperH} with $K'=(2K\mu)^{1/2}$ and by \eqref{comp}, for $v\in E_2(y)$,
\[
\big|\xi\cdot\big(\ba(t,x,v)-\ba(t,x',v)\big)\eta\big|\le C\mu^{1/2}|S_v\xi|\,|S_v\eta|\le C\mu^{1/2}|S_y\xi|\,|S_y\eta| .
\]
Replacing $\xi$ and $\eta$ by $S_y^{-1}\xi$ and $S_y^{-1}\eta$, we get
\begin{equation}\label{cx}
\big|S_y^{-1}\big(\ba(t,x,v)-\ba(t,x',v)\big)S_y^{-1}\big|\le C\varpi(|x-x'|)^{1/2}\qquad(v\in E_2(y)).
\end{equation}

\step{Step 2: dependence on $v$; reduction to a bound for the derivative.}
We claim that for every $y\in\R^3$, every $v\in E_2(y)$ and all $\tau,\xi,\eta\in\R^3$,
\begin{equation}\label{cder}
\big|\xi\cdot(\tau\cdot\nabla_v)\ba(t,x,v)\,\eta\big|\le C\,\big|(r_yS_y)^{-1}\tau\big|\,|S_y\xi|\,|S_y\eta| .
\end{equation}
Suppose this is proved, and let $v,v'\in E_2(y)$. The set $E_2(y)$ is convex, so we may integrate \eqref{cder} with $\tau=v-v'$ along
the segment from $v'$ to $v$. This gives
\begin{equation}\label{cv}
\big|S_y^{-1}\big(\ba(t,x,v)-\ba(t,x,v')\big)S_y^{-1}\big|\le C\big|(r_yS_y)^{-1}(v-v')\big| ,
\end{equation}
and {\rm(c)} follows from \eqref{cx} and \eqref{cv}.

\step{Step 3: proof of \eqref{cder} for $|y|<8$.}
In this case $r_y$ and $S_y$ are comparable to $1$ and to the identity, and $E_2(y)\subset B_{10}$. By {\rm(a)} and \eqref{I1} with
$s=\gamma+1$ we have $|\nabla_v\ba|\le C$ on $B_{10}$, which is \eqref{cder}.

\step{Step 4: proof of \eqref{cder} for $|y|\ge8$.}
Put $e:=\hat y$ and decompose every vector as $\xi=\xi_re+\xi_\perp$ with $\xi_\perp\perp e$. Let $\tilde\tau:=(r_yS_y)^{-1}\tau$. By
\eqref{Sbasic},
\begin{equation}\label{ycomp}
|S_y\xi|^2=\jy^{\gamma}\xi_r^2+\jy^{\gamma+2}|\xi_\perp|^2,\qquad|\tau_\perp|\le|\tilde\tau|,\qquad|\tau|\le|\tilde\tau| .
\end{equation}
Differentiating \eqref{wedge} in the direction $\tau$ gives
\[
\xi\cdot\big((\tau\cdot\nabla)a\big)(z)\,\eta=a_\gamma\Big(\gamma|z|^{\gamma-2}(z\cdot\tau)\,(z\wedge\xi)\cdot(z\wedge\eta)
+|z|^{\gamma}\big[(\tau\wedge\xi)\cdot(z\wedge\eta)+(z\wedge\xi)\cdot(\tau\wedge\eta)\big]\Big).
\]
We estimate the wedge products through their components,
\[
|z\wedge\xi|\le|\xi_r|\,|z_\perp|+|z|\,|\xi_\perp|,\qquad|\tau\wedge\xi|\le|\xi_r|\,|\tau_\perp|+|\tau|\,|\xi_\perp|,
\]
and in the same way for $\eta$. Collecting the terms according to the components of $\xi$ and $\eta$, we obtain
\begin{multline}\label{dera}
\big|\xi\cdot((\tau\cdot\nabla)a)(z)\,\eta\big|\le C|\xi_r|\,|\eta_r|\Big(|z|^{\gamma-1}|z_\perp|^2|\tau|+|z|^{\gamma}|z_\perp|\,|\tau_\perp|\Big)\\
+C|z|^{\gamma+1}|\tau|\,\big(|\xi_r|\,|\eta_\perp|+|\xi_\perp|\,|\eta_r|+|\xi_\perp|\,|\eta_\perp|\big).
\end{multline}
Now let $v\in E_2(y)$ and $z=v-u$. Since $v-y=r_yS_y\tv$ with $|\tv|<2$, we have $|(v-y)_\perp|\le2$, and therefore
\[
|z_\perp|\le|(v-y)_\perp|+|u_\perp|\le3\langle u\rangle,\qquad|z_\perp|\le|z| .
\]
Together with \eqref{ycomp}, this shows that the first bracket in \eqref{dera} is bounded by
\[
2\,|\tilde\tau|\,|z|^{\gamma}|z_\perp|\le6\,|\tilde\tau|\,\min\big\{|z|^{\gamma+1},|z|^{\gamma}\langle u\rangle\big\} .
\]
We multiply \eqref{dera} by $f(u)\le K\langle u\rangle^{-k}$ and integrate in $u$. By \eqref{I2} with $(s,j)=(\gamma,1)$, by \eqref{I1}
with $s=\gamma+1$, and since $\jv$ and $\jy$ are comparable,
\[
\big|\xi\cdot(\tau\cdot\nabla_v)\ba(t,x,v)\,\eta\big|\le C|\tilde\tau|\Big(\jy^{\gamma}|\xi_r|\,|\eta_r|
+\jy^{\gamma+1}\big(|\xi_r|\,|\eta_\perp|+|\xi_\perp|\,|\eta_r|+|\xi_\perp|\,|\eta_\perp|\big)\Big).
\]
By \eqref{ycomp}, each of the four products on the right is at most $|S_y\xi|\,|S_y\eta|$. This proves \eqref{cder}.
\end{proof}

\begin{proof}[Proof of {\rm(d)}]
For $|v|<1$ the claim follows from {\rm(a)}. Let $|v|\ge1$. We bound the components of $\beta$ parallel and orthogonal to $v$
separately.

\step{Step 1: the component parallel to $v$.}
By \eqref{I1} with $s=\gamma+1$, as in the proof of {\rm(a)},
\[
|P_v\beta|\le|\beta|\le C\jv^{\gamma+1} .
\]

\step{Step 2: the component orthogonal to $v$.}
Let $e'$ be a unit vector with $e'\perp v$, and let $z=v-u$. Since $b(z)=-2a_\gamma|z|^{\gamma}z$ and
$e'\cdot z=-e'\cdot P_v^\perp u$,
\[
|e'\cdot b(z)|\le C|z|^{\gamma}|P_v^\perp z|\le C\min\big\{|z|^{\gamma+1},|z|^{\gamma}\langle u\rangle\big\} .
\]
Integrating against $f(u)\le K\langle u\rangle^{-k}$ and using \eqref{I2} with $(s,j)=(\gamma,1)$, we get
\[
|P_v^\perp\beta|\le C\jv^{\gamma} .
\]

\step{Step 3: conclusion.}
Since $S_v^{-1}=\jv^{-\gamma/2}P_v+\jv^{-1-\gamma/2}P_v^\perp$,
\[
\big|r_vS_v^{-1}\beta\big|\le\jv^{-\kappa}\Big(\jv^{-\gamma/2}|P_v\beta|+\jv^{-1-\gamma/2}|P_v^\perp\beta|\Big)
\le C\jv^{-\kappa}\Big(\jv^{1+\gamma/2}+\jv^{-1+\gamma/2}\Big)\le C,
\]
because $1+\gamma/2\le\kappa$.
\end{proof}

\section{The dual identity}\label{sec:dual}
\rev{This section contains two computations in the velocity variable, both for fixed $(t,x)$. The first is behind the formal identity \eqref{formal}. It moves the dependence of the coefficients $\ba^h$ and $\beta^h$ on the density $h$ onto a convolution operator acting on the test function; applied with $h=w$, it turns the contribution of $Q(w,g)$ into the pairing of $w$ with $\K_g\phi$. The second is the integration by parts which shows that every weak solution is a distributional solution. At the end of the section we define the local adjoint operator $\cP$.}

\begin{lemma}\label{lem:dual}
Let $\ell>5+\gamma$, let $h,u:\R^3\to\R$ be measurable with $|h|+|u|\le C_*\jv^{-\ell}$, and let $\phi\in H^2(\R^3)$ have compact support.
Then $\ba^h\in C^1(\R^3)$, $\beta^h\in C(\R^3)$ and $\partial_i\ba^h_{ij}=\beta^h_j$. The following identities hold.
\begin{itemize}
\item[\rm(a)] With
\[
\K_u\phi:=a_{ij}*_v(u\,\partial_{ij}\phi)-2\,b_i*_v(u\,\partial_i\phi)
\]
we have
\begin{equation}\label{d2}
\int_{\R^3}u\,\big(\ba^h:D_v^2\phi+2\beta^h\cdot\nabla_v\phi\big)\dd v=\int_{\R^3}h\,\K_u\phi\dd v .
\end{equation}
\item[\rm(b)] If in addition $u\in H^1_{\rm loc}(\R^3)$, then
\begin{equation}\label{d1}
\int_{\R^3}\big(\ba^h\nabla_vu-\beta^hu\big)\cdot\nabla_v\phi\dd v=-\int_{\R^3}u\,\big(\ba^h:D_v^2\phi+2\beta^h\cdot\nabla_v\phi\big)\dd v .
\end{equation}
\end{itemize}
\end{lemma}
\begin{proof}
\step{Step 1: regularity of the coefficients.}
The kernels satisfy
\[
|a(z)|\le C|z|^{\gamma+2},\qquad|\nabla a(z)|+|b(z)|\le C|z|^{\gamma+1},
\]
so they are locally integrable. Since $\ell>5+\gamma$, the integrals defining $a*h$, $b*h$ and $(\nabla a)*h$ converge absolutely and
locally uniformly in $v$, and they define continuous functions. Hence $\ba^h\in C^1$ with $\nabla_v\ba^h=(\nabla a)*h$, and
\[
\partial_i\ba^h_{ij}=(\partial_ia_{ij})*h=b_j*h=\beta^h_j .
\]

\step{Step 2: proof of {\rm(a)}.}
Inserting $\ba^h=a*h$ and $\beta^h=b*h$, the left-hand side of \eqref{d2} becomes
\[
\int_{\R^3}\!\int_{\R^3}u(v)\Big(a_{ij}(v-v')\,\partial_{ij}\phi(v)+2\,b_i(v-v')\,\partial_i\phi(v)\Big)h(v')\dd v'\dd v .
\]
This double integral converges absolutely: the function $v\mapsto\int(|a|+|b|)(v-v')\,|h(v')|\dd v'$ is bounded on the support of
$\phi$, and $u\,(|D_v^2\phi|+|\nabla_v\phi|)\in L^1(\R^3)$. By Fubini's theorem we may therefore integrate in $v$ first. Since
$a(-z)=a(z)$ and $b(-z)=-b(z)$,
\begin{align*}
\int_{\R^3}a_{ij}(v-v')\,u(v)\,\partial_{ij}\phi(v)\dd v&=\big(a_{ij}*(u\,\partial_{ij}\phi)\big)(v'),\\
\int_{\R^3}b_i(v-v')\,u(v)\,\partial_i\phi(v)\dd v&=-\big(b_i*(u\,\partial_i\phi)\big)(v') .
\end{align*}
Hence the double integral equals $\int h(v')\,\K_u\phi(v')\dd v'$, which is \eqref{d2}.

\step{Step 3: proof of {\rm(b)}.}
Consider the vector field $X:=\ba^h\nabla_v\phi$. By Step 1 it belongs to $H^1(\R^3)$, it has compact support, and
\[
\nabla_v\cdot X=\ba^h:D_v^2\phi+\beta^h\cdot\nabla_v\phi .
\]
Since $u\in H^1_{\rm loc}(\R^3)$ and $\ba^h$ is symmetric, an integration by parts gives
\[
\int_{\R^3}\ba^h\nabla_vu\cdot\nabla_v\phi\dd v=\int_{\R^3}\nabla_vu\cdot X\dd v=-\int_{\R^3}u\,\nabla_v\cdot X\dd v
=-\int_{\R^3}u\,\big(\ba^h:D_v^2\phi+\beta^h\cdot\nabla_v\phi\big)\dd v .
\]
Subtracting $\int u\,\beta^h\cdot\nabla_v\phi\dd v$ from both sides gives \eqref{d1}.
\end{proof}

Part {\rm(b)} of Lemma~\ref{lem:dual} implies that every weak solution $h$ of \eqref{L} is a distributional solution. Indeed, for
a.e.\ $(t,x)$ the function $h(t,x,\cdot)$ belongs to $H^1_{\rm loc}(\R^3)$, so that \eqref{d1} holds with $u=h(t,x,\cdot)$ and
$\phi=\varphi(t,x,\cdot)$. Integrating in $(t,x)$, we see that \eqref{weakform} implies \eqref{dist}.

\medskip
From now on, $f$ and $g$ are as in Theorem~\ref{thm:main}. With $\ba=\ba^f$ and $\beta=\beta^f$ we define the
\emph{local adjoint operator}
\begin{equation}\label{Pdef}
\cP\phi:=Y\phi+\ba:D_v^2\phi+2\beta\cdot\nabla_v\phi=Y\phi+\nabla_v\cdot(\ba\nabla_v\phi)+\beta\cdot\nabla_v\phi .
\end{equation}
The two expressions agree because $\partial_i\ba_{ij}=\beta_j$.

\section{Transport and weak formulations}\label{sec:tools}
\rev{The solutions $\phi$ of the adjoint problem constructed in Section~\ref{sec:solv} are not smooth: we control only $\phi$, $\nabla_v\phi$, $D_v^2\phi$ and $Y\phi$, and only in weighted $L^2$ norms. This section collects the tools for working in that class. Lemma~\ref{lem:T} gives traces in time and an energy identity along the characteristics of $Y$. Theorem~\ref{thm:DY} is the kinetic $L^2$ estimate of \cite{DY}, from which we obtain both the Hessian bound and the solvability of the adjoint problem. Lemma~\ref{lem:dens} admits compactly supported functions of this class as test functions in \eqref{dist}. Finally, Lemma~\ref{lem:pos} is an energy identity for positive parts, which is used only in Section~\ref{sec:cor}.}

\subsection{Transport}
\begin{lemma}\label{lem:T}
Let $J=(s_1,s_2)$ be a bounded interval and $\Omega\subset\R^3$ an open set. Let $u\in L^2(J\times\R^3\times\Omega)$, and suppose that
$Yu$, understood in the sense of distributions, also belongs to $L^2(J\times\R^3\times\Omega)$. Then the function
\[
U(t,y,v):=u(t,y+tv,v)
\]
belongs to $H^1(J;L^2(\R^3\times\Omega))$, and $\partial_tU(t,y,v)=(Yu)(t,y+tv,v)$. In particular $u$ has a representative in
$C(\bar J;L^2(\R^3\times\Omega))$ with the following properties.
\begin{itemize}
\item[\rm(i)] If $\psi$ is bounded and measurable and $\psi(t,y+tv,v)$ does not depend on $t$, then the function
$t\mapsto\iint\psi\,u^2\dd x\dd v$ is absolutely continuous on $\bar J$, with derivative $2\iint\psi\,u\,Yu\dd x\dd v$.
\item[\rm(ii)] If $u(s_2)=0$, then the extension of $u$ by $0$ to $t\ge s_2$ satisfies
$u,Yu\in L^2((s_1,\infty)\times\R^3\times\Omega)$, where $Yu$ is also extended by $0$.
\end{itemize}
\end{lemma}
\begin{proof}
Let $\Theta(t,y,v):=(t,y+tv,v)$. This map is a diffeomorphism of $J\times\R^3\times\Omega$ onto itself, it preserves the Lebesgue
measure, and for every smooth function $\psi$
\begin{equation}\label{conj}
\partial_t(\psi\circ\Theta)=(Y\psi)\circ\Theta .
\end{equation}

\step{Step 1: $U\in H^1(J;L^2)$.}
Let $\psi\in C_c^\infty(J\times\R^3\times\Omega)$ and put $\tilde\psi:=\psi\circ\Theta^{-1}$, so that $(Y\tilde\psi)\circ\Theta=\partial_t\psi$
by \eqref{conj}. Changing variables with $\Theta$ and using the definition of the distribution $Yu$,
\begin{align*}
\iiint U\,\partial_t\psi\dd y\dd v\dd t&=\iiint u\,Y\tilde\psi\dd x\dd v\dd t\\
&=-\iiint(Yu)\,\tilde\psi\dd x\dd v\dd t=-\iiint\big((Yu)\circ\Theta\big)\,\psi\dd y\dd v\dd t .
\end{align*}
Thus $\partial_tU=(Yu)\circ\Theta\in L^2$ in the sense of distributions, that is, $U\in H^1(J;L^2(\R^3\times\Omega))$.

\step{Step 2: continuity of $u$.}
Every element of $H^1(J;L^2)$ has a representative in $C(\bar J;L^2)$. For fixed $t$ we have $u(t)=U(t)\circ\Theta_t^{-1}$, where
$\Theta_t(y,v):=(y+tv,v)$ is a measure preserving shear of $\R^3\times\Omega$. Since the shears act on $L^2(\R^3\times\Omega)$ as a
strongly continuous group of isometries, $u\in C(\bar J;L^2(\R^3\times\Omega))$.

\step{Step 3: proof of {\rm(i)}.}
Let $\psi_0(y,v):=\psi(t,y+tv,v)$, which by assumption does not depend on $t$. Then
\[
\iint\psi\,u^2\dd x\dd v=\iint\psi_0\,U^2\dd y\dd v .
\]
Since $U\in H^1(J;L^2)$ and $\psi_0$ is bounded, the right-hand side is absolutely continuous in $t$ with derivative
$2\iint\psi_0\,U\,\partial_tU\dd y\dd v$. Changing variables back, this derivative equals $2\iint\psi\,u\,Yu\dd x\dd v$.

\step{Step 4: proof of {\rm(ii)}.}
If $u(s_2)=0$, then $U(s_2)=0$. Hence the extension of $U$ by $0$ to $t\ge s_2$ belongs to $H^1((s_1,\infty);L^2)$, and its derivative
is the extension of $\partial_tU$ by $0$. Changing variables back gives {\rm(ii)}.
\end{proof}

For an interval $I\subset(0,T)$ we denote by $\mathcal S(I)$ the set of measurable functions
$\phi$ on $I\times\R^6$ such that
\[
\phi,\ \nabla_v\phi,\ D_v^2\phi,\ Y\phi\in L^2(I\times\R^3\times\Omega)\qquad\text{for every bounded open set }\Omega\subset\R^3_v .
\]
By Lemma~\ref{lem:T}, the traces $\phi(t,\cdot,\cdot)\in L^2_{\rm loc}$ are well defined for every $t\in\bar I$.

\subsection{A kinetic Sobolev estimate}
We state the unweighted $L^2$ estimate of Dong and Yastrzhembskiy with their convention $\partial_t-v\cdot D_x$ for transport.
For $\mathcal T\in(-\infty,\infty]$ let $\R^7_{\mathcal T}:=(-\infty,\mathcal T)\times\R^6$ and
\[
S_2(\R^7_{\mathcal T}):=\big\{u:\ u,\ \partial_tu-v\cdot D_xu,\ D_vu,\ D_v^2u\in L^2(\R^7_{\mathcal T})\big\}.
\]
\begin{theorem}[{\hskip0.0001pt\cite[Theorem~2.6]{DY}}]\label{thm:DY}
Let $\delta\in(0,1)$, $L>0$, and let $\varpi_0:[0,\infty)\to[0,\infty)$ be nondecreasing with $\varpi_0(0+)=0$. Let $a=(a^{ij})$ be
symmetric and measurable and $b=(b^i)$ measurable on $\R^7_{\mathcal T}$. Assume that
\begin{itemize}
\item[\rm(i)] $\delta|\xi|^2\le a^{ij}(z)\xi_i\xi_j\le\delta^{-1}|\xi|^2$ for all $z\in\R^7_{\mathcal T}$ and $\xi\in\R^3$,
\item[\rm(ii)] $|b(z)|\le L$ for all $z\in\R^7_{\mathcal T}$,
\item[\rm(iii)] $|a(t,x_1,v_1)-a(t,x_2,v_2)|\le\varpi_0(|x_1-x_2|+|v_1-v_2|)$ for all $t$, $x_1,x_2$, $v_1,v_2$.
\end{itemize}
Let
\[
Pu:=\partial_tu-v\cdot D_xu-a^{ij}D_{v_iv_j}u .
\]
Then the following two assertions hold.
\begin{itemize}
\item[\rm(1)] There is a constant $N=N(\delta,L,\varpi_0)$, independent of $\mathcal T$, such that for all $u\in S_2(\R^7_{\mathcal T})$,
\[
\|D_v^2u\|_{L^2(\R^7_{\mathcal T})}\le N\Big(\|Pu+b\cdot D_vu\|_{L^2(\R^7_{\mathcal T})}+\|u\|_{L^2(\R^7_{\mathcal T})}\Big).
\]
\item[\rm(2)] If $\mathcal T<\infty$, then for every $S'<\mathcal T$ and every $F\in L^2((S',\mathcal T)\times\R^6)$ there is a unique
function $u$ with $Pu+b\cdot D_vu=F$ in $(S',\mathcal T)\times\R^6$ whose extension by zero to $t\le S'$ belongs to
$S_2(\R^7_{\mathcal T})$.
\end{itemize}
\end{theorem}

\begin{proof}
This is a special case of \cite[Theorem~2.6]{DY}. Condition {\rm(iii)} implies Assumption~2.3~$(\gamma_0)$ of \cite{DY} for every
$\gamma_0>0$, with $R_0=R_0(\varpi_0,\gamma_0)$. Indeed, on a time slice of a kinetic cylinder of radius $r$ one has
\[
|x_1-x_2|\le2r^3,\qquad|v_1-v_2|\le2r,
\]
so the mean oscillation of $a$ over the slice is at most $\varpi_0(2r^3+2r)$. Statement {\rm(1)} follows from the estimate of
\cite[Theorem~2.6(i)]{DY} with $\lambda=\lambda_0$ and the triangle inequality; see also \cite[Corollary~2.8]{DY}. The constants $N$
and $\lambda_0$ in \cite[Theorem~2.6(i)]{DY} do not depend on $\mathcal T$, so the constant in {\rm(1)} is independent of
$\mathcal T$. Statement {\rm(2)} is the existence and uniqueness assertion of \cite[Theorem~2.6(ii)]{DY}, where the zero initial
condition is understood in the sense of \cite[Definition~2.5]{DY}.
\end{proof}

\subsection{Approximation and the positive part}
\begin{lemma}[Density of smooth test functions]\label{lem:dens}
Let $\Phi\in L^2((0,T)\times\R^6)$ vanish outside a compact set $\mathcal C\subset(0,T)\times\R^6$, and suppose that
\[
\nabla_v\Phi,\quad D_v^2\Phi,\quad Y\Phi\quad\text{belong to }L^2((0,T)\times\R^6).
\]
Then there are functions $\Phi_j\in C_c^\infty((0,T)\times\R^6)$, supported in a fixed compact subset of $(0,T)\times\R^6$, such that
\[
\Phi_j\to\Phi,\qquad\nabla_v\Phi_j\to\nabla_v\Phi,\qquad D_v^2\Phi_j\to D_v^2\Phi,\qquad Y\Phi_j\to Y\Phi\qquad\text{in }L^2 .
\]
Consequently, if $h$ is a distributional solution of \eqref{L} in $(0,T)\times\R^6$, then \eqref{dist} holds with $\varphi=\Phi$.
\end{lemma}
\begin{proof}
Let $\varrho_\delta$ and $\varrho'_{\delta'}$ be standard mollifiers on $\R^4_{t,x}$ and
on $\R^3_v$, supported in balls of radii $\delta$ and $\delta'$.

\step{Step 1: regularization in $(t,x)$.}
Put $\Phi^\delta:=\Phi*_{t,x}\varrho_\delta$. Convolution in $(t,x)$ commutes with $\nabla_v$, with $D_v^2$ and with $Y$, because the
coefficients of $Y$ depend only on $v$. Hence, as $\delta\to0$,
\[
\Phi^\delta\to\Phi,\qquad\nabla_v\Phi^\delta\to\nabla_v\Phi,\qquad D_v^2\Phi^\delta\to D_v^2\Phi,\qquad Y\Phi^\delta\to Y\Phi\qquad\text{in }L^2 .
\]
Moreover $\nabla_x\Phi^\delta\in L^2$ for every fixed $\delta$.

\step{Step 2: regularization in $v$.}
Put $\Phi^{\delta,\delta'}:=\Phi^\delta*_v\varrho'_{\delta'}$, which belongs to $C_c^\infty$. As $\delta'\to0$, the functions
$\Phi^{\delta,\delta'}$, $\nabla_v\Phi^{\delta,\delta'}$ and $D_v^2\Phi^{\delta,\delta'}$ converge in $L^2$ to $\Phi^\delta$,
$\nabla_v\Phi^\delta$ and $D_v^2\Phi^\delta$. The operator $Y$ does not commute with the convolution in $v$, but the commutator is
explicit:
\[
Y\Phi^{\delta,\delta'}(t,x,v)-\big((Y\Phi^\delta)*_v\varrho'_{\delta'}\big)(t,x,v)
=\int_{\R^3}(v-u)\cdot\nabla_x\Phi^\delta(t,x,u)\,\varrho'_{\delta'}(v-u)\dd u .
\]
Since $|v-u|\le\delta'$ on the support of $\varrho'_{\delta'}(v-u)$, the $L^2$ norm of the right-hand side is at most
$\delta'\,\|\nabla_x\Phi^\delta\|_{L^2}$. Hence also $Y\Phi^{\delta,\delta'}\to Y\Phi^\delta$ in $L^2$ as $\delta'\to0$.

\step{Step 3: conclusion.}
For small $\delta$ and $\delta'$ all supports lie in a fixed compact subset of $(0,T)\times\R^6$. Choosing first $\delta=\delta_j\to0$
and then $\delta'=\delta'_j$ small enough, we obtain the sequence $\Phi_j:=\Phi^{\delta_j,\delta'_j}$. On the fixed compact set the
functions $h$, $\ba^h$ and $\beta^h$ are bounded. Therefore the left-hand side of \eqref{dist} is continuous with respect to the
convergence just described, and \eqref{dist} for $\varphi=\Phi_j$ implies \eqref{dist} for $\varphi=\Phi$.
\end{proof}

\begin{lemma}[Energy identity for the positive part]\label{lem:pos}
Let $0<s_0<T$ and $\rho>0$. Let
\[
u,\ F\in L^2_{\rm loc}((0,T)\times\R^6),\qquad\nabla_vu,\ G\in L^2_{\rm loc}((0,T)\times\R^6;\R^3)
\]
satisfy
\[
Yu=\nabla_v\cdot G+F\qquad\text{in the sense of distributions in }(0,T)\times\R^6 ,
\]
and suppose that $u\le0$ a.e.\ outside $[s_0,T)\times B_\rho\times B_\rho$. Then $u_+:=\max\{u,0\}$ satisfies
$\nabla_vu_+=\mathbf 1_{\{u>0\}}\nabla_vu$ and, for a.e.\ $t\in(0,T)$,
\begin{equation}\label{pos}
\frac12\iint_{\R^6}u_+(t,x,v)^2\dd x\dd v=\int_0^t\!\!\iint_{\R^6}\big(-G\cdot\nabla_vu_++F\,u_+\big)\dd x\dd v\dd s .
\end{equation}
\end{lemma}
\begin{proof}
The formula for $\nabla_vu_+$ is the standard one, applied for a.e.\ $(t,x)$. It remains to prove \eqref{pos}.

\step{Step 1: regularization.}
Let $\varrho_\delta$ be a mollifier on $\R^4_{t,x}$ supported in $B_\delta$, with $0<\delta<\min\{s_0/2,1\}$. For a function $q$ on
$(0,T)\times\R^6$ we write $q^\delta:=q*_{t,x}\varrho_\delta$, which is defined on $(\delta,T-\delta)\times\R^6$. Since $Y$ and
$\nabla_v$ commute with this convolution,
\[
Yu^\delta=\nabla_v\cdot G^\delta+F^\delta\qquad\text{in }(\delta,T-\delta)\times\R^6 .
\]
The functions $u^\delta$, $\partial_tu^\delta$, $\nabla_xu^\delta$ and $\nabla_vu^\delta$ belong to $L^2_{\rm loc}$. Consequently:
\begin{itemize}
\item $u^\delta\in H^1_{\rm loc}$, and $\nabla_v\cdot G^\delta=Yu^\delta-F^\delta\in L^2_{\rm loc}$;
\item $u^\delta_+\in H^1_{\rm loc}$, with $\nabla_{t,x,v}u^\delta_+=\mathbf 1_{\{u^\delta>0\}}\nabla_{t,x,v}u^\delta$;
\item $u^\delta\le0$ a.e.\ outside $[s_0-\delta,T)\times B_{\rho+1}\times B_\rho$.
\end{itemize}

\step{Step 2: the identity for $u^\delta$.}
By Step 1 the function $t\mapsto\frac12\iint(u^\delta_+)^2\dd x\dd v$ is locally absolutely continuous on $(\delta,T-\delta)$ and
vanishes for $t<s_0-\delta$. Its derivative is
\begin{align*}
\iint u^\delta_+\,\partial_tu^\delta\dd x\dd v&=\iint u^\delta_+\,Yu^\delta\dd x\dd v\\
&=\iint u^\delta_+\big(\nabla_v\cdot G^\delta+F^\delta\big)\dd x\dd v\\
&=\iint\big(-G^\delta\cdot\nabla_vu^\delta_++F^\delta u^\delta_+\big)\dd x\dd v .
\end{align*}
In the first line we used that $\iint v\cdot\nabla_x(u^\delta_+)^2\dd x\dd v=0$, and in the last line we integrated by parts in $v$.
Integrating in time, we see that \eqref{pos} holds for $(u^\delta,G^\delta,F^\delta)$ and every $t\in(\delta,T-\delta)$.

\step{Step 3: the limit $\delta\to0$.}
On compact subsets of $(0,T)\times\R^6$ we have
\[
u^\delta\to u,\qquad\nabla_vu^\delta\to\nabla_vu,\qquad G^\delta\to G,\qquad F^\delta\to F\qquad\text{in }L^2 .
\]
Hence $u^\delta_+\to u_+$ in $L^2$, and $\nabla_vu^\delta_+\to\nabla_vu_+$ weakly in $L^2$, because this sequence is bounded in $L^2$.
All integrands in \eqref{pos} are supported in the compact set $[s_0/2,t]\times\bar B_{\rho+1}\times\bar B_\rho$. Therefore the
right-hand sides converge for every $t\in(0,T)$, as integrals of products of a strongly and a weakly convergent sequence, and the
left-hand sides converge for a.e.\ $t$ along a subsequence. This proves \eqref{pos}.
\end{proof}

\section{The local estimate in Landau-adapted boxes}\label{sec:local}
In Sections \ref{sec:local}--\ref{sec:solv} the assumptions and the notation of Section~\ref{sec:coef} are in force. For matrices,
$|\cdot|$ denotes the Hilbert--Schmidt norm.

\subsection{Cutoff functions adapted to the boxes}
Fix $\zeta\in C_c^\infty(B_1)$ with $0\le\zeta\le1$ and $\zeta=1$ on $B_{1/2}$. For $y\in\R^3$ put
\[
\zeta_y(v):=\zeta\big((r_yS_y)^{-1}(v-y)\big).
\]
Then $\operatorname{supp}\zeta_y\subset E_1(y)$ and $\zeta_y=1$ on $E_{1/2}(y)$. By the chain rule,
\begin{equation}\label{zetader}
S_y\nabla_v\zeta_y=r_y^{-1}\,(\nabla\zeta)\big((r_yS_y)^{-1}(v-y)\big),\qquad
S_yD_v^2\zeta_y\,S_y=r_y^{-2}\,(D^2\zeta)\big((r_yS_y)^{-1}(v-y)\big),
\end{equation}
and in particular
\begin{equation}\label{zetabd}
|S_y\nabla_v\zeta_y|\le Cr_y^{-1},\qquad|S_yD_v^2\zeta_y\,S_y|\le Cr_y^{-2}.
\end{equation}

\subsection{The local estimate}
\begin{proposition}\label{prop:local}
There is a constant $N_0$ such that the following holds for every $y\in\R^3$ and every interval $I=(t_1,t_2)\subset(0,T)$. Suppose that
\[
\phi,\ \nabla_v\phi,\ D_v^2\phi,\ Y\phi,\ F\ \in L^2(I\times\R^3\times E_2(y)),
\]
that $\cP\phi=F$ a.e.\ in $I\times\R^3\times E_2(y)$, and that $\phi(t_2)=0$. Then
\begin{equation}\label{local}
\big\|S_yD_v^2(\zeta_y\phi)S_y\big\|_{L^2(I\times\R^6)}\le N_0\Big(\|F\|+r_y^{-2}\|\phi\|+r_y^{-1}\|S_y\nabla_v\phi\|\Big),
\end{equation}
where the three norms on the right-hand side are taken in $L^2(I\times\R^3\times E_1(y))$.
\end{proposition}

\begin{proof}
Throughout the proof we write
\[
S:=S_y,\qquad r:=r_y,\qquad h:=t_2-t_1,
\]
and $\|\cdot\|$ without a subscript denotes the norm of $L^2(I\times\R^3\times E_1(y))$.

\step{Step 1: the equation for $u:=\zeta_y\phi$.}
We extend $u$ by $0$ for $v\notin E_1(y)$. Then $u,\nabla_vu,D_v^2u,Yu\in L^2(I\times\R^6)$ and $u(t_2)=0$. Since $Y\zeta_y=0$, the
product rule gives
\begin{equation}\label{equ}
\cP u=F_y:=\zeta_yF+\phi\,\big(\ba:D_v^2\zeta_y+2\beta\cdot\nabla_v\zeta_y\big)+2\,\ba\nabla_v\zeta_y\cdot\nabla_v\phi\qquad\text{in }I\times\R^6 .
\end{equation}
To estimate $F_y$ we express the coefficients in the frame $S$. By Lemma~\ref{lem:coef}(b),(d) and \eqref{comp}, for $v\in E_2(y)$,
\begin{equation}\label{normcoef}
c\,\mathrm{Id}\le S^{-1}\ba(t,x,v)S^{-1}\le C\,\mathrm{Id},\qquad |rS^{-1}\beta(t,x,v)|\le C .
\end{equation}
We write the three coefficient terms in \eqref{equ} as
\begin{align*}
\ba:D_v^2\zeta_y&=\tr\big((S^{-1}\ba S^{-1})\,(SD_v^2\zeta_y\,S)\big),\\
\beta\cdot\nabla_v\zeta_y&=(S^{-1}\beta)\cdot(S\nabla_v\zeta_y),\\
\ba\nabla_v\zeta_y\cdot\nabla_v\phi&=(S^{-1}\ba S^{-1})\,S\nabla_v\zeta_y\cdot S\nabla_v\phi .
\end{align*}
By \eqref{zetabd} and \eqref{normcoef},
\[
|\ba:D_v^2\zeta_y|\le Cr^{-2},\qquad|\beta\cdot\nabla_v\zeta_y|\le Cr^{-2},\qquad|\ba\nabla_v\zeta_y\cdot\nabla_v\phi|\le Cr^{-1}|S\nabla_v\phi| .
\]
Since $\zeta_y$ is supported in $E_1(y)$, we conclude that
\begin{equation}\label{Fy}
\|F_y\|_{L^2(I\times\R^6)}\le\|F\|+Cr^{-2}\|\phi\|+Cr^{-1}\|S\nabla_v\phi\| .
\end{equation}

\step{Step 2: rescaling.}
Consider the affine map
\[
\Lambda(\ttil,\tx,\tv):=\big(t_2+r^2\ttil,\ r^3S\tx+r^2\ttil\,y,\ y+rS\tv\big).
\]
It is a bijection of $(-h/r^2,0)\times\R^6$ onto $I\times\R^6$ with constant Jacobian, and it maps the set $\{|\tv|<\rho\}$ onto the
set $\{v\in E_\rho(y)\}$. Let $\tilde u:=u\circ\Lambda$. By the chain rule,
\begin{equation}\label{chain}
\nabla_{\tv}\tilde u=r\,(S\nabla_vu)\circ\Lambda,\qquad
D^2_{\tv}\tilde u=r^2\,(SD_v^2u\,S)\circ\Lambda,\qquad
(\partial_{\ttil}+\tv\cdot\nabla_{\tx})\tilde u=r^2\,(Yu)\circ\Lambda .
\end{equation}
The third identity follows by adding
\[
\partial_{\ttil}\tilde u=r^2\,(\partial_tu+y\cdot\nabla_xu)\circ\Lambda\qquad\text{and}\qquad
\tv\cdot\nabla_{\tx}\tilde u=r^3\,(S\tv\cdot\nabla_xu)\circ\Lambda=r^2\,\big((v-y)\cdot\nabla_xu\big)\circ\Lambda .
\]
We multiply \eqref{equ} by $r^2$, compose with $\Lambda$ and use \eqref{chain}. This gives
\[
(\partial_{\ttil}+\tv\cdot\nabla_{\tx})\tilde u+\bar A:D^2_{\tv}\tilde u+\bar B\cdot\nabla_{\tv}\tilde u=r^2F_y\circ\Lambda ,
\]
with the rescaled coefficients
\[
\bar A:=(S^{-1}\ba S^{-1})\circ\Lambda,\qquad\bar B:=2r\,(S^{-1}\beta)\circ\Lambda .
\]
The bounds \eqref{normcoef} hold for $|\tv|<2$. Since $\tilde u$ vanishes for $|\tv|\ge1$, we may modify the coefficients
for $|\tv|\ge1$. Let $\chi\in C_c^\infty(B_2)$ with $0\le\chi\le1$ and $\chi=1$ on $B_1$, and put
\[
a^\sharp:=\chi(\tv)\,\bar A+(1-\chi(\tv))\,\mathrm{Id},\qquad b^\sharp:=\chi(\tv)\,\bar B .
\]
Then
\begin{equation}\label{eqsharp}
(\partial_{\ttil}+\tv\cdot\nabla_{\tx})\tilde u+a^\sharp:D^2_{\tv}\tilde u+b^\sharp\cdot\nabla_{\tv}\tilde u=r^2F_y\circ\Lambda
\qquad\text{in }(-h/r^2,0)\times\R^6 .
\end{equation}

\step{Step 3: the hypotheses of Theorem~\ref{thm:DY} hold uniformly in $y$ and $I$.}
We check conditions {\rm(i)--(iii)} of Theorem~\ref{thm:DY} for the coefficients $a^\sharp$ and $b^\sharp$.

\emph{Ellipticity and boundedness.} By \eqref{normcoef} there are $\delta\in(0,1)$ and $L>0$, independent of $y$ and $I$, such that
\[
\delta\,\mathrm{Id}\le a^\sharp\le\delta^{-1}\mathrm{Id},\qquad|b^\sharp|\le L .
\]

\emph{Continuity in $(\tx,\tv)$.} Fix $\ttil$, let $|\tv_1|,|\tv_2|<2$ and $\tx_1,\tx_2\in\R^3$, and let
$(t,x_i,v_i):=\Lambda(\ttil,\tx_i,\tv_i)$. \rev{These two points have the same time $t$ and the same Galilean shift $r^2\ttil\,y$, which cancels in $x_1-x_2$. Hence}
\[
v_i\in E_2(y),\qquad\big|(rS)^{-1}(v_1-v_2)\big|=|\tv_1-\tv_2|,\qquad|x_1-x_2|=r^3|S(\tx_1-\tx_2)|\le|\tx_1-\tx_2| .
\]
By Lemma~\ref{lem:coef}(c),
\[
\big|\bar A(\ttil,\tx_1,\tv_1)-\bar A(\ttil,\tx_2,\tv_2)\big|\le C\Big(|\tv_1-\tv_2|+\varpi(|\tx_1-\tx_2|)^{1/2}\Big).
\]
Now $a^\sharp-\mathrm{Id}=\chi\,(\bar A-\mathrm{Id})$, where $\chi$ is Lipschitz and supported in $B_2$, and $\bar A$ is bounded on
$\{|\tv|<2\}$. Hence $a^\sharp$ satisfies condition {\rm(iii)} of Theorem~\ref{thm:DY} with the modulus
\[
\varpi_0(\rho):=C\rho+C\varpi(\rho)^{1/2},
\]
which does not depend on $y$, $I$ or $\ttil$.

\step{Step 4: application of Theorem~\ref{thm:DY}.}
Theorem~\ref{thm:DY} is stated for the forward operator, so we reverse time. Let
\[
V(s,\tx,\tv):=\tilde u(-s,\tx,\tv),\qquad s\in(0,h/r^2).
\]
By \eqref{eqsharp},
\[
\partial_sV-\tv\cdot D_{\tx}V-a^\sharp:D^2_{\tv}V-b^\sharp\cdot D_{\tv}V=-r^2\,(F_y\circ\Lambda)(-s,\cdot),\qquad V(0)=0,
\]
where the coefficients are evaluated at time $-s$. We extend $V$ and the right-hand side by $0$, and $(a^\sharp,b^\sharp)$ by
$(\mathrm{Id},0)$, to $s\le0$. By Lemma~\ref{lem:T}(ii) (after reversing time), the extended function $V$ belongs to
$S_2(\R^7_{h/r^2})$ and satisfies the equation in $\R^7_{h/r^2}$; the extended coefficients still satisfy {\rm(i)--(iii)}. Therefore
Theorem~\ref{thm:DY}(1) yields
\[
\|D^2_{\tv}V\|_{L^2}\le N\Big(\|r^2F_y\circ\Lambda\|_{L^2}+\|V\|_{L^2}\Big),\qquad N=N(\delta,L,\varpi_0).
\]
We change variables back with $\Lambda$. The Jacobian is the same constant in all three terms, and by \eqref{chain} we get
\[
r^2\,\big\|SD_v^2u\,S\big\|_{L^2(I\times\R^6)}\le N\Big(r^2\,\|F_y\|_{L^2(I\times\R^6)}+\|u\|_{L^2(I\times\R^6)}\Big).
\]
Dividing by $r^2$ and using $\|u\|_{L^2(I\times\R^6)}\le\|\phi\|$ together with \eqref{Fy}, we arrive at \eqref{local}.
\end{proof}

\section{Energy estimate and solvability of the local adjoint problem}\label{sec:solv}

\subsection{Truncated operators}
Let $\eta\in C_c^\infty(B_2)$ with $0\le\eta\le1$ and $\eta=1$ on $B_1$, and for $R\ge1$
let $\eta_R(v):=\eta(v/R)$. We define
\[
A^R:=\eta_R\,\ba+(1-\eta_R)\,\mathrm{Id},\qquad b^R:=\eta_R\,\beta,\qquad
\cP^R\phi:=Y\phi+\nabla_v\cdot(A^R\nabla_v\phi)+b^R\cdot\nabla_v\phi .
\]
By \eqref{div} and \eqref{Pdef}, $\cP^R\phi=\cP\phi$ on $\{|v|<R\}$. In nondivergence form,
\begin{equation}\label{PRnd}
\cP^R\phi=Y\phi+A^R:D_v^2\phi+\tilde b^R\cdot\nabla_v\phi,\qquad
\tilde b^R_j:=2\eta_R\,\beta_j+\partial_i\eta_R\,(\ba_{ij}-\delta_{ij}).
\end{equation}

\subsection{The weighted energy estimate}
\begin{lemma}[Weighted energy estimate]\label{lem:energy}
There are constants $h_0>0$ and $C$, independent of $R$, such that the following holds. Let $I=(t_1,t_2)\subset(0,T)$ with
$t_2-t_1\le h_0$, and let
\[
\phi,\ \nabla_v\phi,\ D_v^2\phi,\ Y\phi\in L^2(I\times\R^6)
\]
satisfy $\cP^R\phi=F$ a.e.\ in $I\times\R^6$ and $\phi(t_2)=0$. Then
\begin{equation}\label{energy}
\sup_{t\in I}\|\omega_{m_1}\phi(t)\|_{L^2(\R^6)}^2+\iiint_{I\times\R^6}\omega_{m_1}^2\,A^R\nabla_v\phi\cdot\nabla_v\phi\dd v\dd x\dd t
\le32\,\|\omega_{m_1}F\|_{L^1(I;L^2(\R^6))}^2 .
\end{equation}
\end{lemma}

\begin{proof}
We work with the bounded weight
\[
\omega:=\omega_{m_1,\epsilon},\qquad\epsilon\in(0,1],
\]
and let $\epsilon\to0$ at the end. By {\rm(W2)}, $\nabla_v\omega$ is bounded as well. We write $\|\cdot\|$ for the norm of $L^2(\R^6)$
and
\[
D(t):=\iint_{\R^6}\omega^2\,A^R\nabla_v\phi\cdot\nabla_v\phi\dd x\dd v\ \ge0 .
\]

\step{Step 1: the energy identity.}
By {\rm(W1)}, the function $\omega^2(t,y+tv,v)$ does not depend on $t$. Hence Lemma~\ref{lem:T}(i) applies with $\psi=\omega^2$: the
function $t\mapsto\|\omega\phi(t)\|^2$ is absolutely continuous on $\bar I$, and
\begin{equation}\label{en1}
\frac{\dd}{\dd t}\|\omega\phi\|^2=2\iint\omega^2\phi\,Y\phi
=2\iint\omega^2\phi\,F-2\iint\omega^2\phi\,\nabla_v\cdot(A^R\nabla_v\phi)-2\iint\omega^2\phi\,b^R\cdot\nabla_v\phi .
\end{equation}

\step{Step 2: the diffusion term.}
We integrate by parts in $v$ and use the Cauchy--Schwarz inequality for the nonnegative matrix $A^R$:
\begin{align*}
-2\iint\omega^2\phi\,\nabla_v\cdot(A^R\nabla_v\phi)
&=2\iint\omega^2A^R\nabla_v\phi\cdot\nabla_v\phi+4\iint\omega\,\phi\,A^R\nabla_v\phi\cdot\nabla_v\omega\\
&\ge\iint\omega^2A^R\nabla_v\phi\cdot\nabla_v\phi-4\iint\phi^2A^R\nabla_v\omega\cdot\nabla_v\omega .
\end{align*}
By Lemma~\ref{lem:coef}(b) and {\rm(W2)},
\[
A^R\nabla_v\omega\cdot\nabla_v\omega\le C_1|S_v\nabla_v\omega|^2+|\nabla_v\omega|^2\le C\omega^2 .
\]
Therefore
\begin{equation}\label{en2}
-2\iint\omega^2\phi\,\nabla_v\cdot(A^R\nabla_v\phi)\ge D(t)-C\|\omega\phi\|^2 .
\end{equation}

\step{Step 3: the drift term.}
Since $2\phi\,b^R\cdot\nabla_v\phi=b^R\cdot\nabla_v(\phi^2)$, an integration by parts gives
\[
-2\iint\omega^2\phi\,b^R\cdot\nabla_v\phi=\iint\nabla_v\cdot(\omega^2b^R)\,\phi^2 .
\]
By \eqref{div} we have $\nabla_v\cdot\beta=-\bc$, and hence
\[
\nabla_v\cdot(\omega^2b^R)=-\omega^2\eta_R\,\bc+\omega^2\beta\cdot\nabla_v\eta_R+2\,\omega\,\eta_R\,\beta\cdot\nabla_v\omega .
\]
Each of the three terms is bounded from below by $-C\omega^2$:
\begin{itemize}
\item for the first term, because $0\le\bc\le C$ by Lemma~\ref{lem:coef}(a);
\item for the second term, because $\nabla_v\eta_R$ is supported in $\{R\le|v|\le2R\}$ and bounded by $C/R$, so that by
Lemma~\ref{lem:coef}(a)
\[
|\beta\cdot\nabla_v\eta_R|\le C\langle2R\rangle^{(\gamma+1)_+}R^{-1}\le C ;
\]
\item for the third term, because by Lemma~\ref{lem:coef}(a) and {\rm(W2)}
\[
|\beta\cdot\nabla_v\omega|\le C\jv^{(\gamma+1)_+}\big(\jv^{-\kappa}+\jv^{-1}\big)\,\omega\le C\omega ,
\]
where we used $(\gamma+1)_+\le\kappa$.
\end{itemize}
Therefore
\begin{equation}\label{en3}
-2\iint\omega^2\phi\,b^R\cdot\nabla_v\phi\ge-C\|\omega\phi\|^2 .
\end{equation}

\step{Step 4: a differential inequality.}
We insert \eqref{en2} and \eqref{en3} into \eqref{en1} and estimate the term with $F$ by the Cauchy--Schwarz inequality. This gives a
constant $C_E$, independent of $R$ and $\epsilon$, such that
\[
\frac{\dd}{\dd t}\|\omega\phi\|^2\ge D(t)-C_E\|\omega\phi\|^2-2\,\|\omega\phi\|\,\|\omega F\| .
\]
Let $M:=\sup_{t\in I}\|\omega\phi(t)\|$, which is finite because $\omega$ is bounded and $\phi\in C(\bar I;L^2(\R^6))$. Integrating
over $(t,t_2)$ and using $\phi(t_2)=0$, we obtain for every $t\in I$
\begin{equation}\label{en4}
\|\omega\phi(t)\|^2+\int_t^{t_2}D(s)\dd s\le C_E\,(t_2-t_1)\,M^2+2M\,\|\omega F\|_{L^1(I;L^2)} .
\end{equation}

\step{Step 5: conclusion.}
Choose $h_0$ such that $C_Eh_0\le\frac12$. Taking the supremum over $t$ in \eqref{en4} gives
$M^2\le\frac12M^2+2M\|\omega F\|_{L^1(I;L^2)}$, that is,
\[
M\le4\,\|\omega F\|_{L^1(I;L^2)} .
\]
Inserting this bound in the right-hand side of \eqref{en4} shows that
\[
\sup_{t\in I}\|\omega\phi(t)\|^2\le16\,\|\omega F\|_{L^1(I;L^2)}^2\qquad\text{and}\qquad\int_ID(s)\dd s\le16\,\|\omega F\|_{L^1(I;L^2)}^2 .
\]
This is \eqref{energy} with $\omega_{m_1,\epsilon}$ in place of $\omega_{m_1}$.
Since $\omega_{m_1,\epsilon}$ increases to $\omega_{m_1}$ as $\epsilon\searrow0$, the monotone convergence theorem gives \eqref{energy}.
\end{proof}

\subsection{Solvability}
\begin{proposition}[Solvability]\label{prop:solv}
There is a constant $C$ such that for every interval $I=(t_1,t_2)\subset(0,T)$ with $t_2-t_1\le h_0$ and every
$\Psi\in C_c^\infty(I\times\R^6)$ there exists $\phi\in\mathcal S(I)$ with
\[
\cP\phi=\Psi\quad\text{a.e.\ in }I\times\R^6,\qquad\phi(t_2)=0,
\]
which satisfies the estimates
\begin{align}
\sup_{t\in I}\|\omega_{m_1}\phi(t)\|_{L^2(\R^6)}+\big\|\omega_{m_1}|S_v\nabla_v\phi|\big\|_{L^2(I\times\R^6)}&\le C\,\|\omega_{m_1}\Psi\|_{L^1(I;L^2(\R^6))},\label{E1}\\
\big\|\omega_{m}|S_vD_v^2\phi\,S_v|\big\|_{L^2(I\times\R^6)}&\le C\,\|\omega_{m_1}\Psi\|_{L^2(I\times\R^6)} .\label{E2}
\end{align}
\end{proposition}

\begin{proof}
\step{Step 1: the truncated problems.}
Fix $R\ge1$. We verify the hypotheses of Theorem~\ref{thm:DY} for the operator $\cP^R$ in the nondivergence form \eqref{PRnd}.
\begin{itemize}
\item The matrix $A^R$ is bounded and uniformly elliptic, with constants depending on $R$. This follows from Lemma~\ref{lem:coef}(b),
because $S_v$ is comparable to the identity for $|v|\le2R$.
\item The drift $\tilde b^R$ is bounded, by Lemma~\ref{lem:coef}(a).
\item The matrix $A^R$ satisfies condition {\rm(iii)} with some modulus $\varpi_0^R$. Indeed, $A^R=\mathrm{Id}$ for $|v|\ge2R$,
while for $|v|\le2R$ Lemma~\ref{lem:coef} gives
\[
|\nabla_v\ba|\le C_R,\qquad|\ba(t,x,v)-\ba(t,x',v)|\le C_R\,\varpi(|x-x'|)^{1/2}.
\]
\end{itemize}
After the change of variable $s=t_2-t$, the problem $\cP^R\phi=\Psi$, $\phi(t_2)=0$ becomes a Cauchy problem with zero initial
condition for an operator of the form considered in Theorem~\ref{thm:DY} (the coefficients are extended by $(\mathrm{Id},0)$ to $s\le0$,
as in Step~4 of the proof of Proposition~\ref{prop:local}). By Theorem~\ref{thm:DY}(2) there exists a function $\phi^R$ with
\[
\phi^R,\ \nabla_v\phi^R,\ D_v^2\phi^R,\ Y\phi^R\in L^2(I\times\R^6),\qquad\cP^R\phi^R=\Psi\ \text{a.e.},\qquad\phi^R(t_2)=0 .
\]

\step{Step 2: bounds which are uniform in $R$.}
By Lemma~\ref{lem:energy} with $F=\Psi$,
\begin{equation}\label{enR}
\sup_{t\in I}\|\omega_{m_1}\phi^R(t)\|_{L^2(\R^6)}^2+\iiint_{I\times\R^6}\omega_{m_1}^2\,A^R\nabla_v\phi^R\cdot\nabla_v\phi^R\dd v\dd x\dd t
\le32\,\|\omega_{m_1}\Psi\|_{L^1(I;L^2(\R^6))}^2 .
\end{equation}
We now derive local bounds for the second derivatives. Fix $R_0\ge1$ and let $R\ge R_0+4$. For $|y|\le R_0+1$ we have
$E_2(y)\subset B_{R_0+3}$, and on $B_{R_0+3}$
\[
\cP^R=\cP,\qquad A^R=\ba\ge c_1S_v^2,\qquad\omega_{m_1}\ge c(R_0)>0 .
\]
Hence \eqref{enR} bounds $\phi^R$ and $S_v\nabla_v\phi^R$ in $L^2(I\times\R^3\times B_{R_0+3})$ by a constant $C(R_0,\Psi)$. We apply
Proposition~\ref{prop:local} to $\phi^R$ in the box $E_2(y)$, with $F=\Psi$. Since $\zeta_y=1$ on $E_{1/2}(y)$, and since $r_y$ and
$S_y$ are bounded from above and from below for $|y|\le R_0+1$, the estimate \eqref{local} gives
\[
\|D_v^2\phi^R\|_{L^2(I\times\R^3\times E_{1/2}(y))}\le C(R_0,\Psi)\qquad(|y|\le R_0+1).
\]
The ball $B_{R_0}$ is covered by finitely many sets $E_{1/2}(y_i)$ with $|y_i|\le R_0$. Consequently
\[
\phi^R,\quad\nabla_v\phi^R,\quad D_v^2\phi^R,\quad Y\phi^R=\Psi-\ba:D_v^2\phi^R-2\beta\cdot\nabla_v\phi^R
\]
are bounded in $L^2(I\times\R^3\times B_{R_0})$, uniformly in $R\ge R_0+4$.

\step{Step 3: passage to the limit.}
By Step 2 and a diagonal argument there are $R_j\to\infty$ and $\phi\in\mathcal S(I)$ such that, for every $R_0\ge1$,
\[
\phi^{R_j}\rightharpoonup\phi,\qquad\nabla_v\phi^{R_j}\rightharpoonup\nabla_v\phi,\qquad D_v^2\phi^{R_j}\rightharpoonup D_v^2\phi,\qquad
Y\phi^{R_j}\rightharpoonup Y\phi
\]
weakly in $L^2(I\times\R^3\times B_{R_0})$.
We check the three properties of $\phi$.

\emph{The equation.} Since $\cP^{R_j}=\cP$ on $\{|v|<R_j\}$ and the coefficients of $\cP$ are bounded on $B_{R_0}$, we can pass to
the limit in $\cP\phi^{R_j}=\Psi$ on $I\times\R^3\times B_{R_0}$. Thus $\cP\phi=\Psi$ a.e.\ in $I\times\R^6$.

\emph{The terminal condition.} In the variables of Lemma~\ref{lem:T}, the functions $\Phi^{R_j}(t,y,v):=\phi^{R_j}(t,y+tv,v)$
converge weakly in $H^1(I;L^2(\R^3\times B_{R_0}))$ to $\Phi(t,y,v):=\phi(t,y+tv,v)$. Therefore
$\Phi^{R_j}(t)\rightharpoonup\Phi(t)$ weakly in $L^2(\R^3\times B_{R_0})$ for every $t\in\bar I$, and in particular $\phi(t_2)=0$.

\emph{The energy bound.} In the variables $(t,y,v)$ the weight $\omega_{m_1}$ becomes
$\langle y\jv^{-\kappa}\rangle^{n_x}\jv^{-m_1}$, which does not depend on $t$. By weak lower semicontinuity (applied first with the
weight and the domain truncated, and then combined with monotone convergence), the bound \eqref{enR} passes to the limit, with
$A^R$ replaced by $\ba$:
\[
\sup_{t\in I}\|\omega_{m_1}\phi(t)\|_{L^2(\R^6)}^2+\iiint_{I\times\R^6}\omega_{m_1}^2\,\ba\nabla_v\phi\cdot\nabla_v\phi\dd v\dd x\dd t
\le32\,\|\omega_{m_1}\Psi\|_{L^1(I;L^2(\R^6))}^2 .
\]
Since $\ba\ge c_1S_v^2$ by Lemma~\ref{lem:coef}(b), this is \eqref{E1}.

\step{Step 4: the equation for $\te\phi$.}
It remains to prove \eqref{E2}. We apply Proposition~\ref{prop:local} in every box $E_2(y)$ to the function
\[
\phi^\epsilon:=\te\,\phi,\qquad\epsilon\in(0,1],
\]
and then sum over the boxes with the help of \eqref{aver}. All constants below are independent of $\epsilon$.

Since $\te$ and its derivatives in $v$ are bounded and $Y\te=0$, we have $\phi^\epsilon\in\mathcal S(I)$ and $\phi^\epsilon(t_2)=0$.
By the product rule,
\[
\cP\phi^\epsilon=\te\Psi+G^\epsilon,\qquad
G^\epsilon:=\phi\,\big(\ba:D_v^2\te+2\beta\cdot\nabla_v\te\big)+2\,\ba\nabla_v\te\cdot\nabla_v\phi .
\]
By {\rm(W2)} and Lemma~\ref{lem:coef}(a),(b),
\begin{equation}\label{Geps}
|G^\epsilon|\le C\te\,\big(|\phi|+|S_v\nabla_v\phi|\big).
\end{equation}

\step{Step 5: the estimate in one box.}
Fix $y\in\R^3$. On $E_1(y)$, by the product rule, {\rm(W2)}, \eqref{zetabd}, \eqref{comp}
and $r_y\le1$,
\begin{equation}\label{prodrule}
\big|S_yD_v^2(\zeta_y\phi^\epsilon)S_y-\zeta_y\te\,S_yD_v^2\phi\,S_y\big|+r_y^{-1}|S_y\nabla_v\phi^\epsilon|
\le C\te\,\Big(r_y^{-1}|S_y\nabla_v\phi|+r_y^{-2}|\phi|\Big).
\end{equation}
We apply Proposition~\ref{prop:local} to $\phi^\epsilon$, with $F=\te\Psi+G^\epsilon$. Using \eqref{Geps}, \eqref{prodrule},
\eqref{comp} and the fact that $\zeta_y=1$ on $E_{1/2}(y)$, we obtain from \eqref{local}
\begin{equation}\label{onebox}
\big\|\te\,|S_vD_v^2\phi\,S_v|\big\|_{L^2(I\times\R^3\times E_{1/2}(y))}^2
\le C\,\Big\|\te\,\big(|\Psi|+r_v^{-2}|\phi|+r_v^{-1}|S_v\nabla_v\phi|\big)\Big\|_{L^2(I\times\R^3\times E_{1}(y))}^2 .
\end{equation}

\step{Step 6: summation over the boxes.}
We multiply \eqref{onebox} by $\jy^{-2m}$, which is comparable to $\jv^{-2m}$ on $E_1(y)$,
and integrate with respect to the measure $\det(r_yS_y)^{-1}\dd y$. Applying \eqref{aver} for fixed $(t,x)$ to both sides, we get
\[
\big\|\omega_{m,\epsilon}\,|S_vD_v^2\phi\,S_v|\big\|_{L^2(I\times\R^6)}
\le C\,\Big\|\omega_{m,\epsilon}\,\big(|\Psi|+r_v^{-2}|\phi|+r_v^{-1}|S_v\nabla_v\phi|\big)\Big\|_{L^2(I\times\R^6)} .
\]
Since $\omega_{m,\epsilon}\,r_v^{-2}=\omega_{m_1,\epsilon}\le\omega_{m_1}$ and $r_v\le1$, this implies
\begin{multline}\label{summed}
\big\|\omega_{m,\epsilon}\,|S_vD_v^2\phi\,S_v|\big\|_{L^2(I\times\R^6)}\\
\le C\,\Big(\|\omega_{m_1}\Psi\|_{L^2(I\times\R^6)}+\|\omega_{m_1}\phi\|_{L^2(I\times\R^6)}+\big\|\omega_{m_1}|S_v\nabla_v\phi|\big\|_{L^2(I\times\R^6)}\Big).
\end{multline}

\step{Step 7: conclusion.}
By \eqref{E1} and the Cauchy--Schwarz inequality in $t$ (recall that $t_2-t_1\le h_0$),
\[
\|\omega_{m_1}\phi\|_{L^2(I\times\R^6)}+\big\|\omega_{m_1}|S_v\nabla_v\phi|\big\|_{L^2(I\times\R^6)}
\le C\,\|\omega_{m_1}\Psi\|_{L^1(I;L^2(\R^6))}\le C\,\|\omega_{m_1}\Psi\|_{L^2(I\times\R^6)} .
\]
Hence the right-hand side of \eqref{summed} is bounded by $C\|\omega_{m_1}\Psi\|_{L^2(I\times\R^6)}$. Letting $\epsilon\to0$ and using
the monotone convergence theorem, we obtain \eqref{E2}.
\end{proof}

\section{The nonlocal term}\label{sec:nonlocal}
In this section we estimate the operator
\[
\K_g\phi=a_{ij}*_v(g\,\partial_{ij}\phi)-2\,b_i*_v(g\,\partial_i\phi)
\]
for fixed $(t,x)$. We write $L^2_v:=L^2(\R^3_v)$.

\begin{lemma}[Kernel bound]\label{lem:ker}
Let $\mathsf k=\mathsf k(z)$ be a measurable function on $\R^3$ with
\[
|\mathsf k(z)|\le C_{\mathsf k}\,|z|^{-\sigma}\langle z\rangle^{\sigma+2\kappa}\qquad\text{for some }0\le\sigma<3 .
\]
Then there is a constant $C$, depending also on $C_{\mathsf k}$ and $\sigma$, such that the following holds for all $t\in[0,T]$ and
$x\in\R^3$. If $G=G(v)$ is measurable and $\theta\jv^{n}G\in L^2_v$, then the integral defining $\mathsf k*_vG(v)$ converges absolutely
for a.e.\ $v\in\R^3$, and
\begin{equation}\label{ker}
\|\omega_{m_1}\,\mathsf k*_vG\|_{L^2_v}\le C\,\|\theta\jv^{n}G\|_{L^2_v}.
\end{equation}
This applies to the kernels
\begin{align*}
\mathsf k&=a_{ij}&&\text{with }\sigma=(-\gamma-2)_+,\\
\mathsf k&=b_i&&\text{with }\sigma=(-\gamma-1)_+,\\
\mathsf k&=|z|^{\gamma}&&\text{with }\sigma=-\gamma\quad(\text{if }\gamma>-3).
\end{align*}
\end{lemma}
\begin{proof}
Fix $(t,x)$ and write $\theta(v)$ for $\theta(t,x,v)$. By {\rm(W3)} and the assumption on $\mathsf k$, for all $v,u\in\R^3$,
\begin{equation}\label{ker1}
\omega_{m_1}(v)\,|\mathsf k(v-u)|\le C\,\theta(u)\,\jv^{-m_1}\,|v-u|^{-\sigma}\,\langle v-u\rangle^{\sigma+2\kappa+n_x} .
\end{equation}
We split the convolution into the contributions of $\{|v-u|<1\}$ and $\{|v-u|\ge1\}$,
\[
\mathsf k*_vG=T_1+T_2,\qquad T_1(v):=\int_{|v-u|<1}\mathsf k(v-u)\,G(u)\dd u,\qquad T_2(v):=\int_{|v-u|\ge1}\mathsf k(v-u)\,G(u)\dd u .
\]

\step{Step 1: the near part $T_1$.}
For $|v-u|<1$ the right-hand side of \eqref{ker1} is at most $C\theta(u)\,|v-u|^{-\sigma}$. Hence
\[
\omega_{m_1}(v)\int_{|v-u|<1}|\mathsf k(v-u)|\,|G(u)|\dd u\le C\int_{|v-u|<1}|v-u|^{-\sigma}\,\theta(u)\,|G(u)|\dd u .
\]
Since $|z|^{-\sigma}\mathbf 1_{\{|z|<1\}}$ is integrable on $\R^3$, Young's inequality shows that the right-hand side belongs to
$L^2_v$. In particular, the integral defining $T_1(v)$ converges absolutely for a.e.\ $v$, and
\[
\|\omega_{m_1}T_1\|_{L^2_v}\le C\,\|\theta G\|_{L^2_v}\le C\,\|\theta\jv^{n}G\|_{L^2_v} .
\]

\step{Step 2: the far part $T_2$.}
Let $\varsigma:=2\kappa+n_x<6$. For $|v-u|\ge1$ we have
\[
|v-u|^{-\sigma}\langle v-u\rangle^{\sigma+\varsigma}\le C\langle v-u\rangle^{\varsigma}\le C\jv^{\varsigma}\langle u\rangle^{\varsigma},
\]
so that, by \eqref{ker1} and the Cauchy--Schwarz inequality,
\begin{align*}
\omega_{m_1}(v)\int_{|v-u|\ge1}|\mathsf k(v-u)|\,|G(u)|\dd u
&\le C\jv^{\varsigma-m_1}\int_{\R^3}\theta(u)\langle u\rangle^{\varsigma}|G(u)|\dd u\\
&\le C\jv^{\varsigma-m_1}\,\big\|\langle u\rangle^{\varsigma-n}\big\|_{L^2}\,\|\theta\jv^{n}G\|_{L^2_v} .
\end{align*}
The functions $\langle u\rangle^{\varsigma-n}$ and $\jv^{\varsigma-m_1}$ belong to $L^2(\R^3)$, because $n-\varsigma>\frac32$ and
$m_1-\varsigma>\frac32$. Hence the integral defining $T_2(v)$ converges absolutely for every $v$, and taking the $L^2$ norm in $v$
gives $\|\omega_{m_1}T_2\|_{L^2_v}\le C\|\theta\jv^{n}G\|_{L^2_v}$. Together with Step~1, this proves \eqref{ker}.

\step{Step 3: the kernels.}
We have $|a_{ij}(z)|\le C|z|^{\gamma+2}$ and $|b_i(z)|\le C|z|^{\gamma+1}$. If the exponent $\gamma+2$ (respectively $\gamma+1$) is
nonnegative, it is at most $2\kappa$ and we may take $\sigma=0$. If it is negative, we take $\sigma$ equal to its absolute value.
For $\mathsf k=|z|^{\gamma}$ with $\gamma>-3$ we take $\sigma=-\gamma<3$.
\end{proof}

\begin{lemma}[The nonlocal term]\label{lem:K}
Let $|g|\le K\jv^{-k}$ with $k\ge44$, and let $\varpi_g$ be a modulus of continuity of $g$ in $v$ which is uniform in $(t,x)$. For
every $\delta>0$ there is a constant $C_\delta$ such that for all $(t,x)$ and all $\phi=\phi(v)\in H^2_{\rm loc}(\R^3)$ for which the
right-hand side is finite,
\begin{equation}\label{Kbound}
\|\omega_{m_1}\K_g\phi\|_{L^2_v}\le\delta\,\big\|\omega_m|S_vD_v^2\phi\,S_v|\big\|_{L^2_v}+C_\delta\|\omega_{m_1}\phi\|_{L^2_v}+C\big\|\omega_{m_1}|S_v\nabla_v\phi|\big\|_{L^2_v}.
\end{equation}
\end{lemma}

\begin{proof}
We fix $(t,x)$; all norms are norms of $L^2_v$. Let
\[
\mathcal N(\phi):=\big\|\omega_m|S_vD_v^2\phi\,S_v|\big\|+\|\omega_{m_1}\phi\|+\big\|\omega_{m_1}|S_v\nabla_v\phi|\big\|
\]
be the sum of the three norms on the right-hand side of \eqref{Kbound}.

\step{Step 1: a crude bound.}
By \eqref{Sbasic},
\[
|D_v^2\phi|\le\jv^{-\gamma}|S_vD_v^2\phi\,S_v|,\qquad|\nabla_v\phi|\le\jv^{-\gamma/2}|S_v\nabla_v\phi| .
\]
Since $|g|\le K\jv^{-k}$ and $k\ge n+m_1+3\ge n+m-\gamma$, it follows that
\[
\theta\jv^{n}|g|\,|D_v^2\phi|\le K\,\omega_m|S_vD_v^2\phi\,S_v|,\qquad
\theta\jv^{n}|g|\,|\nabla_v\phi|\le K\,\omega_{m_1}|S_v\nabla_v\phi| .
\]
Hence Lemma~\ref{lem:ker}, applied to the kernels $a_{ij}$ and $b_i$, gives
\begin{equation}\label{crude}
\|\omega_{m_1}\,a_{ij}*(g\,\partial_{ij}\phi)\|+\|\omega_{m_1}\,b_i*(g\,\partial_i\phi)\|
\le C\,\big\|\theta\jv^{n-k}\big(|D_v^2\phi|+|\nabla_v\phi|\big)\big\|\le C\,\mathcal N(\phi).
\end{equation}

\step{Step 2: reduction to smooth $\phi$.}
By \eqref{crude} the map $\phi\mapsto\omega_{m_1}\K_g\phi$ is continuous with respect to $\mathcal N$. We claim that
$C_c^\infty(\R^3)$ is dense, with respect to $\mathcal N$, in the set of all $\phi\in H^2_{\rm loc}(\R^3)$ with
$\mathcal N(\phi)<\infty$. Then it suffices to prove \eqref{Kbound} for $\phi\in C_c^\infty(\R^3)$.

To prove the claim, let $\chi\in C_c^\infty(B_2)$ with $\chi=1$ on $B_1$, and let $\chi_L(v):=\chi(v/L)$ for $L\ge1$. The derivatives
of $\chi_L$ are supported in $\{L\le|v|\le2L\}$, where $|S_v\xi|\le CL^{\kappa}|\xi|$ by \eqref{Sbasic}. Hence
\[
|S_v\nabla_v\chi_L|\le CL^{\kappa-1}\le C,\qquad|S_vD_v^2\chi_L\,S_v|\le CL^{2\kappa-2}\le C .
\]
By the product rule and $\omega_m\le\omega_{m_1}$,
\[
\mathcal N(\phi-\chi_L\phi)\le C\Big(\big\|\mathbf 1_{\{|v|\ge L\}}\,\omega_m|S_vD_v^2\phi\,S_v|\big\|
+\big\|\mathbf 1_{\{|v|\ge L\}}\,\omega_{m_1}\phi\big\|+\big\|\mathbf 1_{\{|v|\ge L\}}\,\omega_{m_1}|S_v\nabla_v\phi|\big\|\Big).
\]
The right-hand side tends to zero as $L\to\infty$, by the dominated convergence theorem. For fixed $L$, the weights $\omega_m$,
$\omega_{m_1}$ and the matrices $S_v$ are bounded on $B_{2L+1}$, so that
\[
\mathcal N(\psi)\le C_L\,\|\psi\|_{H^2(\R^3)}\qquad\text{for every }\psi\in H^2(\R^3)\text{ supported in }B_{2L+1} .
\]
Hence the mollifications of $\chi_L\phi$ converge to $\chi_L\phi$ with respect to $\mathcal N$. This proves the claim.

\step{Step 3: the first-order part.}
By the second half of \eqref{crude},
\[
\|\omega_{m_1}\,b_i*(g\,\partial_i\phi)\|\le C\,\big\|\omega_{m_1}|S_v\nabla_v\phi|\big\| ,
\]
which is the last term of \eqref{Kbound}. It remains to estimate the second-order part $a_{ij}*(g\,\partial_{ij}\phi)$.

\step{Step 4: splitting of $g$.}
Let $\epsilon\in(0,1]$ and let $g_\epsilon:=g*_v\varrho_\epsilon$, where $\varrho_\epsilon$ is a standard mollifier on $\R^3_v$
supported in $B_\epsilon$. Then
\[
|g-g_\epsilon|\le\varpi_g(\epsilon)\qquad\text{and}\qquad|g-g_\epsilon|\le CK\jv^{-k} .
\]
Taking the geometric mean of these two bounds and using $k\ge2(n+m-\gamma)$, we get
\begin{equation}\label{gsmall}
\jv^{n+m-\gamma}\,|g-g_\epsilon|\le C\,\varpi_g(\epsilon)^{1/2} .
\end{equation}
Moreover, the mollified function satisfies
\begin{equation}\label{gsmooth}
|D_v^jg_\epsilon|\le C\epsilon^{-j}\jv^{-k}\qquad(j=0,1,2).
\end{equation}
We write
\[
a_{ij}*(g\,\partial_{ij}\phi)=a_{ij}*\big((g-g_\epsilon)\,\partial_{ij}\phi\big)+a_{ij}*(g_\epsilon\,\partial_{ij}\phi)
\]
and estimate the two terms separately.

\step{Step 5: the rough part.}
By Lemma~\ref{lem:ker}, by $|D_v^2\phi|\le\jv^{-\gamma}|S_vD_v^2\phi\,S_v|$ and by \eqref{gsmall},
\begin{align*}
\big\|\omega_{m_1}\,a_{ij}*\big((g-g_\epsilon)\,\partial_{ij}\phi\big)\big\|
&\le C\,\big\|\theta\jv^{n}\,(g-g_\epsilon)\,D_v^2\phi\big\|\\
&\le C\,\big\|\jv^{n+m-\gamma}|g-g_\epsilon|\cdot\omega_m|S_vD_v^2\phi\,S_v|\big\|\\
&\le C\,\varpi_g(\epsilon)^{1/2}\,\big\|\omega_m|S_vD_v^2\phi\,S_v|\big\| .
\end{align*}

\step{Step 6: the smooth part.}
By the product rule,
\[
g_\epsilon\,\partial_{ij}\phi=\partial_{ij}(g_\epsilon\phi)-\partial_i(\partial_jg_\epsilon\,\phi)-\partial_j(\partial_ig_\epsilon\,\phi)+\partial_{ij}g_\epsilon\,\phi .
\]
In the sense of distributions, $\partial_ia_{ij}=b_j$ and
\[
\partial_jb_j=-c_\gamma|z|^{\gamma}\quad\text{if }\gamma>-3,\qquad\partial_jb_j=-c_{-3}\,\delta_0\quad\text{if }\gamma=-3 .
\]
Convolving the previous identity with $a_{ij}$ and moving the derivatives to the kernel, we obtain
\[
a_{ij}*(g_\epsilon\,\partial_{ij}\phi)=(\partial_jb_j)*(g_\epsilon\phi)-2\,b_j*(\partial_jg_\epsilon\,\phi)+a_{ij}*(\partial_{ij}g_\epsilon\,\phi).
\]
Each of the three terms contains $\phi$ without derivatives. By Lemma~\ref{lem:ker} for the kernels $|z|^{\gamma}$, $b_j$ and $a_{ij}$
(and directly for the term $c_{-3}\,g_\epsilon\phi$ when $\gamma=-3$), by \eqref{gsmooth} and by $k\ge n+m_1$,
\[
\big\|\omega_{m_1}\,a_{ij}*(g_\epsilon\,\partial_{ij}\phi)\big\|\le C\epsilon^{-2}\,\big\|\theta\jv^{n-k}\phi\big\|\le C\epsilon^{-2}\,\|\omega_{m_1}\phi\| .
\]

\step{Step 7: conclusion.}
Collecting Steps 3, 5 and 6, we have shown that for every $\epsilon\in(0,1]$
\[
\|\omega_{m_1}\K_g\phi\|\le C\,\varpi_g(\epsilon)^{1/2}\,\big\|\omega_m|S_vD_v^2\phi\,S_v|\big\|+C\epsilon^{-2}\,\|\omega_{m_1}\phi\|
+C\,\big\|\omega_{m_1}|S_v\nabla_v\phi|\big\| .
\]
Given $\delta>0$, we choose $\epsilon$ so small that $C\varpi_g(\epsilon)^{1/2}\le\delta$ and put $C_\delta:=C\epsilon^{-2}$. This is
\eqref{Kbound}.
\end{proof}

\section{Testing the difference of the solutions}\label{sec:weak}
\rev{We now make the formal identity \eqref{formal} rigorous for functions $\phi$ with the bounds of Proposition~\ref{prop:solv}. Such functions are neither smooth nor compactly supported, and they are controlled only in the weighted norms of \eqref{E1} and \eqref{E2}. The pairing with $w$ therefore puts the weight $\omega_{m_1}$ on $\phi$ and the dual weight $\omega_{m_1}^{-1}$ on $w$. We begin by checking that $w$ is square integrable against the dual weight, although it need not decay in $x$.}

Let $f$ and $g$ be as in Theorem~\ref{thm:main}, and let $w:=f-g$, so that $|w|\le2K\jv^{-k}$. By {\rm(W3)} we have
$\langle x\rangle^{4}\le C\jv^{8}\theta$, and therefore
\[
\omega_{m_1}^{-1}=\theta^{-1}\jv^{8}\le C\langle x\rangle^{-4}\jv^{16},\qquad\text{so that}\qquad
\omega_{m_1}^{-1}\jv^{4}\,|w|\le C\,\langle x\rangle^{-4}\jv^{-2}\cdot\jv^{22}|w| .
\]
Since $\langle x\rangle^{-4}\jv^{-2}\in L^2(\R^6)$, and since $\jv^{22}|w|\le(2K)^{1/2}\jv^{22-k/2}|w|^{1/2}\le C|w|^{1/2}$ for
$k\ge44$, we obtain for every $t\in[0,T]$
\begin{equation}\label{wdual}
\big\|\omega_{m_1}^{-1}\jv^{4}\,w(t)\big\|_{L^2(\R^6)}\le C\,\|\jv^{22}w(t)\|_{L^\infty(\R^6)}\le C\,\|w(t)\|_{L^\infty(\R^6)}^{1/2} .
\end{equation}

\begin{lemma}\label{lem:weak}
Let $I=(t_1,t_2)\subset(0,T)$, and let $\phi\in\mathcal S(I)$ satisfy $\omega_{m_1}\cP\phi\in L^2(I\times\R^6)$ and
\begin{equation}\label{phiclass}
\sup_{t\in I}\|\omega_{m_1}\phi(t)\|_{L^2(\R^6)}+\big\|\omega_{m_1}|S_v\nabla_v\phi|\big\|_{L^2(I\times\R^6)}+\big\|\omega_{m}|S_vD_v^2\phi\,S_v|\big\|_{L^2(I\times\R^6)}<\infty .
\end{equation}
Then, for every $\tau\in(t_1,t_2)$,
\begin{equation}\label{weak}
\iint_{\R^6}w\,\phi\dd x\dd v\Big|_{t=\tau}^{t=t_2}=\int_\tau^{t_2}\!\!\iint_{\R^6}w\,\big(\cP\phi+\K_g\phi\big)\dd x\dd v\dd t ,
\end{equation}
and all integrals converge absolutely.
\end{lemma}

\begin{proof}
\step{Step 1: a cutoff which commutes with $Y$.}
Let $\chi\in C_c^\infty(B_2)$ with $0\le\chi\le1$ and $\chi=1$ on $B_1$. For $L\ge1$ put
\[
\chi_L(t,x,v):=\chi\Big(\frac{x-tv}{L}\Big)\,\chi\Big(\frac{v}{L}\Big),\qquad\phi_L:=\chi_L\,\phi .
\]
Then
\[
Y\chi_L=0,\qquad|\nabla_v\chi_L|\le\frac{C}{L},\qquad|D_v^2\chi_L|\le\frac{C}{L^2},
\]
and $\phi_L$, $\nabla_v\phi_L$, $D_v^2\phi_L$, $Y\phi_L$ belong to $L^2(I\times\R^6)$. In the variables of Lemma~\ref{lem:T}, the
function
\[
\Phi_L(t,y,v):=\phi_L(t,y+tv,v)=\chi\Big(\frac yL\Big)\,\chi\Big(\frac vL\Big)\,\phi(t,y+tv,v)
\]
belongs to $H^1(I;L^2)$ and vanishes outside a fixed compact set of $(y,v)$. Since $f$ and $g$ are bounded and uniformly continuous
on $[0,T]\times\R^6$, it follows that for $h\in\{f,g\}$ the function
\begin{equation}\label{tcont}
t\longmapsto\iint_{\R^6}h\,\phi_L\dd x\dd v=\iint_{\R^6}h(t,y+tv,v)\,\Phi_L(t,y,v)\dd y\dd v
\end{equation}
is continuous on $\bar I$.

\step{Step 2: a time-sliced identity for $f$ and for $g$.}
Let $h\in\{f,g\}$. We claim that for every $\tau\in(t_1,t_2)$
\begin{equation}\label{vw}
\iint_{\R^6}h\,\phi_L\dd x\dd v\Big|_{t=\tau}^{t=t_2}
=\int_\tau^{t_2}\!\!\iint_{\R^6}h\,\big(Y\phi_L+\ba^h:D_v^2\phi_L+2\beta^h\cdot\nabla_v\phi_L\big)\dd x\dd v\dd t .
\end{equation}
To prove this, let $t_1<\tau<\tau'<t_2$ and $0<\epsilon<(\tau'-\tau)/2$, and let $\varsigma_\epsilon=\varsigma_\epsilon(t)$ be the
piecewise affine function which vanishes outside $(\tau,\tau')$, equals $1$ on $[\tau+\epsilon,\tau'-\epsilon]$, and is affine on
$[\tau,\tau+\epsilon]$ and on $[\tau'-\epsilon,\tau']$. The function $\varsigma_\epsilon\phi_L$, extended by $0$, vanishes outside a
compact subset of $(0,T)\times\R^6$, and
\[
Y(\varsigma_\epsilon\phi_L)=\varsigma_\epsilon\,Y\phi_L+\varsigma_\epsilon'\,\phi_L .
\]
By Lemma~\ref{lem:dens}, the distributional formulation \eqref{dist} holds with $\varphi=\varsigma_\epsilon\phi_L$. Since
$\varsigma_\epsilon'=1/\epsilon$ on $(\tau,\tau+\epsilon)$ and $\varsigma_\epsilon'=-1/\epsilon$ on $(\tau'-\epsilon,\tau')$, this
reads
\begin{multline*}
\iiint h\,\varsigma_\epsilon\,\big(Y\phi_L+\ba^h:D_v^2\phi_L+2\beta^h\cdot\nabla_v\phi_L\big)\dd v\dd x\dd t\\
+\frac1\epsilon\int_\tau^{\tau+\epsilon}\!\!\iint h\,\phi_L\dd x\dd v\dd t-\frac1\epsilon\int_{\tau'-\epsilon}^{\tau'}\!\!\iint h\,\phi_L\dd x\dd v\dd t=0 .
\end{multline*}
We let $\epsilon\to0$. By the continuity of the function \eqref{tcont}, the last two terms converge to the values of
$\iint h\,\phi_L\dd x\dd v$ at $t=\tau$ and at $t=\tau'$. Hence \eqref{vw} holds with $\tau'$ in place of $t_2$. Finally, both sides
are continuous in $\tau'$ up to $\tau'=t_2$ (the integrand on the right-hand side is integrable on $I\times\R^6$), and \eqref{vw}
follows.

\step{Step 3: the identity for $w$ and $\phi_L$.}
We subtract \eqref{vw} for $h=g$ from \eqref{vw} for $h=f$. Since $h\mapsto\ba^h$ and $h\mapsto\beta^h$ are linear,
\[
f\,\ba^f-g\,\ba^g=w\,\ba^f+g\,\ba^w,\qquad f\,\beta^f-g\,\beta^g=w\,\beta^f+g\,\beta^w .
\]
Recalling the definition \eqref{Pdef} of $\cP$, we obtain
\[
\iint w\,\phi_L\dd x\dd v\Big|_{t=\tau}^{t=t_2}
=\int_\tau^{t_2}\!\!\iint w\,\cP\phi_L\dd x\dd v\dd t
+\int_\tau^{t_2}\!\!\iint g\,\big(\ba^w:D_v^2\phi_L+2\beta^w\cdot\nabla_v\phi_L\big)\dd x\dd v\dd t .
\]
In the last term we apply Lemma~\ref{lem:dual}(a) with $(h,u)=(w,g)$ for a.e.\ $(t,x)$:
\[
\int_{\R^3}g\,\big(\ba^w:D_v^2\phi_L+2\beta^w\cdot\nabla_v\phi_L\big)\dd v=\int_{\R^3}w\,\K_g\phi_L\dd v .
\]
This gives
\begin{equation}\label{weakL}
\iint w\,\phi_L\dd x\dd v\Big|_{t=\tau}^{t=t_2}=\int_\tau^{t_2}\!\!\iint w\,\big(\cP\phi_L+\K_g\phi_L\big)\dd x\dd v\dd t .
\end{equation}

\step{Step 4: the limit $L\to\infty$ in the local term.}
Since $Y\chi_L=0$, the product rule and Lemma~\ref{lem:coef}(a) give
\[
|\cP\phi_L-\chi_L\,\cP\phi|\le\frac{C}{L}\,\jv^{2\kappa}\big(|\phi|+|\nabla_v\phi|\big)\le\frac{C}{L}\,\jv^{4}\big(|\phi|+|S_v\nabla_v\phi|\big).
\]
By \eqref{wdual}, \eqref{phiclass} and the Cauchy--Schwarz inequality,
\[
\int_\tau^{t_2}\!\!\iint|w|\,|\cP\phi_L-\chi_L\,\cP\phi|\dd x\dd v\dd t\longrightarrow0 .
\]
Moreover, $w\,\cP\phi$ is integrable on $I\times\R^6$, by \eqref{wdual} and the assumption $\omega_{m_1}\cP\phi\in L^2(I\times\R^6)$.
Hence $\int_\tau^{t_2}\!\iint w\,\chi_L\,\cP\phi\to\int_\tau^{t_2}\!\iint w\,\cP\phi$ by the dominated convergence theorem. Combining
the two limits, we get
\[
\int_\tau^{t_2}\!\!\iint w\,\cP\phi_L\dd x\dd v\dd t\longrightarrow\int_\tau^{t_2}\!\!\iint w\,\cP\phi\dd x\dd v\dd t .
\]

\step{Step 5: the limit $L\to\infty$ in the nonlocal term.}
We have
\[
|D_v^2(\phi_L-\phi)|+|\nabla_v(\phi_L-\phi)|\le(1-\chi_L)\big(|D_v^2\phi|+|\nabla_v\phi|\big)+\frac{C}{L}\big(|\nabla_v\phi|+|\phi|\big).
\]
By the first inequality in \eqref{crude}, applied to $\phi_L-\phi$ for a.e.\ $(t,x)$,
\[
\big\|\omega_{m_1}\K_g(\phi_L-\phi)\big\|_{L^2(I\times\R^6)}
\le C\,\Big\|\theta\jv^{n-k}\big(|D_v^2(\phi_L-\phi)|+|\nabla_v(\phi_L-\phi)|\big)\Big\|_{L^2(I\times\R^6)} .
\]
As in Step~1 of the proof of Lemma~\ref{lem:K},
\[
\theta\jv^{n-k}\big(|D_v^2\phi|+|\nabla_v\phi|+|\phi|\big)\le\omega_m|S_vD_v^2\phi\,S_v|+\omega_{m_1}|S_v\nabla_v\phi|+\omega_{m_1}|\phi| ,
\]
and the right-hand side belongs to $L^2(I\times\R^6)$ by \eqref{phiclass}. Hence
$\omega_{m_1}\K_g\phi_L\to\omega_{m_1}\K_g\phi$ in $L^2(I\times\R^6)$, by the dominated convergence theorem. Together with
\eqref{wdual}, this gives
\[
\int_\tau^{t_2}\!\!\iint w\,\K_g\phi_L\dd x\dd v\dd t\longrightarrow\int_\tau^{t_2}\!\!\iint w\,\K_g\phi\dd x\dd v\dd t .
\]

\step{Step 6: the boundary terms and conclusion.}
For $t=\tau$ and $t=t_2$ we have
\[
\iint w(t)\,\phi_L(t)\dd x\dd v\longrightarrow\iint w(t)\,\phi(t)\dd x\dd v ,
\]
by \eqref{wdual} and dominated convergence. \rev{At $t=t_2$ the domination uses the bound}
\[
\|\omega_{m_1}\phi(t_2)\|_{L^2(\R^6)}\le\sup_{t\in I}\|\omega_{m_1}\phi(t)\|_{L^2(\R^6)} ,
\]
which follows from Fatou's lemma, because $\phi(t)\to\phi(t_2)$ in $L^2_{\rm loc}$ as $t\to t_2$.

Letting $L\to\infty$ in \eqref{weakL}, we obtain \eqref{weak}. The absolute convergence of the integrals in \eqref{weak} follows from
\eqref{wdual}, \eqref{phiclass}, \eqref{crude} and the assumption $\omega_{m_1}\cP\phi\in L^2(I\times\R^6)$.
\end{proof}

\section{Proof of Theorem~\ref{thm:main}}\label{sec:proof}
Let $f$ and $g$ satisfy {\rm(A1)--(A2)} with $k\ge44$, and let $w=f-g$. All results of Sections \ref{sec:coef}--\ref{sec:weak}
apply, where $\varpi$ and $\varpi_g$ are moduli of continuity of $f$ and $g$. For an interval $I\subset(0,T)$ we use the weighted space
\[
\cZ(I):=L^2\big(I\times\R^6;\omega_{m_1}^2\dd x\dd v\dd t\big),\qquad\|\Psi\|_{\cZ(I)}:=\|\omega_{m_1}\Psi\|_{L^2(I\times\R^6)} .
\]
Its dual space $\cZ(I)'$ with respect to the $L^2$ pairing is $L^2(I\times\R^6;\omega_{m_1}^{-2}\dd x\dd v\dd t)$, and $w$ belongs
to it: by \eqref{wdual},
\begin{equation}\label{wfin}
\|w\|_{\cZ(I)'}=\|\omega_{m_1}^{-1}w\|_{L^2(I\times\R^6)}<\infty .
\end{equation}

\begin{proof}[Proof of Theorem~\ref{thm:main}]
\step{Step 1: the nonlocal term is a contraction.}
We claim that there is $h_1\in(0,h_0]$ with the following property. For every interval $I=(t_1,t_2)\subset(0,T)$ with
$t_2-t_1\le h_1$ and every $\Psi\in C_c^\infty(I\times\R^6)$, the function $\phi$ of Proposition~\ref{prop:solv} satisfies
\begin{equation}\label{contr}
\|\K_g\phi\|_{\cZ(I)}\le\frac12\,\|\Psi\|_{\cZ(I)} .
\end{equation}
To prove the claim, let $h:=t_2-t_1$. By \eqref{E1}, \eqref{E2} and the Cauchy--Schwarz inequality in $t$,
\begin{align*}
\|\omega_{m_1}\phi\|_{L^2(I\times\R^6)}&\le h^{1/2}\sup_{t\in I}\|\omega_{m_1}\phi(t)\|_{L^2(\R^6)}\le Ch\,\|\Psi\|_{\cZ(I)},\\
\big\|\omega_{m_1}|S_v\nabla_v\phi|\big\|_{L^2(I\times\R^6)}&\le C\,\|\omega_{m_1}\Psi\|_{L^1(I;L^2(\R^6))}\le Ch^{1/2}\,\|\Psi\|_{\cZ(I)},\\
\big\|\omega_{m}|S_vD_v^2\phi\,S_v|\big\|_{L^2(I\times\R^6)}&\le C\,\|\Psi\|_{\cZ(I)} .
\end{align*}
Let $\delta>0$. We take the $L^2$ norm of \eqref{Kbound} with respect to $(t,x)\in I\times\R^3$ and insert the three bounds above:
\[
\|\K_g\phi\|_{\cZ(I)}\le\big(C\delta+C_\delta\,Ch+Ch^{1/2}\big)\,\|\Psi\|_{\cZ(I)} .
\]
We first choose $\delta$ such that $C\delta\le\frac16$, and then $h_1$ such that $C_\delta\,Ch_1+Ch_1^{1/2}\le\frac13$. This gives
\eqref{contr}.

\step{Step 2: uniqueness on a short interval.}
Let $I=(t_1,t_2)\subset(0,T)$ with $t_2-t_1\le h_1$, and suppose that $w(t_1,\cdot,\cdot)=0$; here $t_1=0$ is allowed. We show that
$w=0$ on $[t_1,t_2]\times\R^6$.

Let $\Psi\in C_c^\infty(I\times\R^6)$ and let $\phi$ be the function of Proposition~\ref{prop:solv}. Then $\cP\phi=\Psi$,
$\phi(t_2)=0$, and by \eqref{E1} and \eqref{E2} the function $\phi$ satisfies the hypotheses of Lemma~\ref{lem:weak}. Let
$\tau\in(t_1,t_2)$ be so close to $t_1$ that $\Psi=0$ for $t\le\tau$. Then \eqref{weak} reads
\begin{equation}\label{dualid}
-\iint_{\R^6}w(\tau)\,\phi(\tau)\dd x\dd v=\iiint_{I\times\R^6}w\,\Psi\dd x\dd v\dd t+\int_\tau^{t_2}\!\!\iint_{\R^6}w\,\K_g\phi\dd x\dd v\dd t .
\end{equation}
We let $\tau\searrow t_1$. By \eqref{wdual} and \eqref{E1},
\[
\Big|\iint_{\R^6}w(\tau)\,\phi(\tau)\dd x\dd v\Big|\le\|\omega_{m_1}^{-1}w(\tau)\|_{L^2(\R^6)}\,\|\omega_{m_1}\phi(\tau)\|_{L^2(\R^6)}
\le C\,\|w(\tau)\|_{L^\infty}^{1/2}\,\|\omega_{m_1}\Psi\|_{L^1(I;L^2(\R^6))} .
\]
The right-hand side tends to $0$, because $w$ is uniformly continuous on $[0,T]\times\R^6$ and $w(t_1)=0$. Hence \eqref{dualid}
becomes, in the limit,
\[
\iiint_{I\times\R^6}w\,\Psi\dd x\dd v\dd t=-\iiint_{I\times\R^6}w\,\K_g\phi\dd x\dd v\dd t .
\]
By the Cauchy--Schwarz inequality and \eqref{contr},
\[
\Big|\iiint_{I\times\R^6}w\,\Psi\dd x\dd v\dd t\Big|\le\|w\|_{\cZ(I)'}\,\|\K_g\phi\|_{\cZ(I)}
\le\frac12\,\|w\|_{\cZ(I)'}\,\|\Psi\|_{\cZ(I)} .
\]
This holds for every $\Psi\in C_c^\infty(I\times\R^6)$, and $C_c^\infty(I\times\R^6)$ is dense in $\cZ(I)$. Taking the supremum
over all such $\Psi$ with $\|\Psi\|_{\cZ(I)}\le1$, we find
\[
\|w\|_{\cZ(I)'}\le\frac12\,\|w\|_{\cZ(I)'} .
\]
Since $\|w\|_{\cZ(I)'}$ is finite by \eqref{wfin}, we conclude that $w=0$ a.e.\ in $I\times\R^6$, and hence everywhere on
$[t_1,t_2]\times\R^6$ by continuity.

\step{Step 3: conclusion.}
By {\rm(A1)} we have $w(0,\cdot,\cdot)=0$. Applying Step 2 successively on the intervals
\[
\big(jh_1,\ \min\{(j+1)h_1,T\}\big),\qquad j=0,1,2,\dots,
\]
we obtain $w\equiv0$ on $[0,T]\times\R^6$. This concludes the proof.
\end{proof}

\section{Proof of Corollary~\ref{cor}}\label{sec:cor}
We begin with the propagation of the decay of $g$.

\begin{lemma}[Propagation of polynomial decay]\label{lem:decay}
Let $\ell=5+\gamma+\eta$ for some $\eta>0$. Suppose that $g\ge0$ is uniformly continuous on $[0,T]\times\R^6$, belongs to
$L^{\infty,\ell}([0,T]\times\R^6)$, and is a weak solution of \eqref{L} in $(0,T)\times\R^6$. If
$g(0,\cdot,\cdot)\in L^{\infty,k}(\R^6)$, then $g\in L^{\infty,k}([0,T]\times\R^6)$.
\end{lemma}

\begin{proof}
\step{Step 1: the coefficients of $g$.}
Since $g\ge0$, we have $\ba^g\ge0$ and $\bc^g\ge0$. We claim that there is a constant $K'$, depending only on
$\|g\|_{L^{\infty,\ell}}$, $\gamma$ and $\ell$, such that
\begin{equation}\label{gcoef}
|\ba^g|\le K'\jv^{2\kappa},\qquad|\beta^g|\le K'\jv^{(\gamma+1)_+},\qquad0\le\bc^g\le K' .
\end{equation}
Let $-3<s\le2$ with $\ell>3+s$. If $s\ge0$, then $|v-u|^{s}\le C(\jv^{s}+\langle u\rangle^{s})$, and hence
\[
\int_{\R^3}|v-u|^{s}\langle u\rangle^{-\ell}\dd u\le C\jv^{s}\int_{\R^3}\langle u\rangle^{-\ell}\dd u+C\int_{\R^3}\langle u\rangle^{s-\ell}\dd u\le C\jv^{s} .
\]
If $s<0$, then $|z|^{s}$ and $\langle z\rangle^{-\ell}$ are radially symmetric and decreasing. Writing each of them as a
superposition of indicator functions of centered balls, and using $|B_r(v)\cap B_{r'}(0)|\le|B_r(0)\cap B_{r'}(0)|$, we see that their
convolution attains its maximum at the origin:
\[
\int_{\R^3}|v-u|^{s}\langle u\rangle^{-\ell}\dd u\le\int_{\R^3}|u|^{s}\langle u\rangle^{-\ell}\dd u<\infty .
\]
In both cases the integral is at most $C\jv^{s_+}$. Since $\ell>5+\gamma$, we may take $s=\gamma+2$ and $s=\gamma+1$. This gives the
bounds for $\ba^g$ and $\beta^g$ in \eqref{gcoef}. If $\gamma>-3$, we may take $s=\gamma$, which gives the bound for $\bc^g$. If
$\gamma=-3$, then $\bc^g=c_{-3}\,g$ is bounded.

\step{Step 2: the transport operator in divergence form.}
For a smooth function $\zeta$ let
\[
\mathcal L\zeta:=Y\zeta-\ba^g:D_v^2\zeta-\bc^g\zeta .
\]
By \eqref{div}, $\nabla_v\cdot(\ba^g\nabla_v\zeta-\beta^g\zeta)=\ba^g:D_v^2\zeta+\bc^g\zeta$ in the sense of distributions. Hence
\begin{equation}\label{Lzeta}
Y\zeta=\nabla_v\cdot\big(\ba^g\nabla_v\zeta-\beta^g\zeta\big)+\mathcal L\zeta .
\end{equation}

\step{Step 3: two barriers.}
We use the barriers of \cite[Lemma 2.3]{HST}. Let
\[
\bar\phi:=Me^{\beta_0t}\jv^{-k},\qquad M:=\|g(0)\|_{L^{\infty,k}} .
\]
Then $\bar\phi(0)\ge g(0)$ and $Y\bar\phi=\beta_0\bar\phi$. Since $|D_v^2\bar\phi|\le C\jv^{-2}\bar\phi$ and $2\kappa\le2$, the
bounds \eqref{gcoef} give $|\ba^g:D_v^2\bar\phi|\le CK'\bar\phi$, so that
\[
\mathcal L\bar\phi\ge\big(\beta_0-(C+1)K'\big)\,\bar\phi\ge0\qquad\text{if }\beta_0\ge(C+1)K' .
\]
Next, fix $\ell'\in(0,\ell)$ and let
\[
\psi:=e^{\beta_1t}\,\jv^{-\ell'}\,\big\langle(x-tv)\jv^{-\kappa}\big\rangle .
\]
As in the proof of Lemma~\ref{lem:weights}, with the exponent $n_x$ replaced by $1$,
\[
Y\psi=\beta_1\psi,\qquad|D_v^2\psi|\le C\jv^{-2\kappa}\psi .
\]
By \eqref{gcoef}, it follows that $|\ba^g:D_v^2\psi|\le CK'\psi$, and therefore
\[
\mathcal L\psi\ge\big(\beta_1-(C+1)K'\big)\,\psi\ge\psi\qquad\text{if }\beta_1\ge(C+1)K'+1 .
\]

\step{Step 4: the function $u$ and the support of its positive part.}
Fix $\varepsilon>0$ and put
\[
u:=g-\bar\phi-\varepsilon\psi .
\]
We claim that there are $\rho>0$ and $s_0>0$ such that
\begin{equation}\label{usupp}
u\le0\qquad\text{outside }[s_0,T]\times B_\rho\times B_\rho .
\end{equation}
We have $g\le C\jv^{-\ell}$, while $\psi\jv^{\ell}\to\infty$ as $|x|+|v|\to\infty$, uniformly in $t\in[0,T]$. Hence there is
$\rho>0$ such that
\[
u<0\qquad\text{outside }[0,T]\times B_\rho\times B_\rho .
\]
On $\bar B_\rho\times\bar B_\rho$ we have $u(0)\le-\varepsilon\psi(0)\le-2c$ for some $c>0$, because $g(0)\le\bar\phi(0)$. Since $u$
is uniformly continuous on $[0,T]\times\bar B_\rho\times\bar B_\rho$, there is $s_0>0$ such that
\[
u\le-c\qquad\text{on }[0,s_0]\times\bar B_\rho\times\bar B_\rho .
\]
This proves \eqref{usupp}.

\step{Step 5: the equation for $u$.}
By the weak formulation \eqref{weakform} for $g$ and by \eqref{Lzeta} for $\zeta=\bar\phi+\varepsilon\psi$,
\[
Yu=\nabla_v\cdot G+F,\qquad G:=\ba^g\nabla_vu-\beta^gu,\qquad F:=-\mathcal L(\bar\phi+\varepsilon\psi)\le0,
\]
in the sense of distributions in $(0,T)\times\R^6$. Moreover, $u$, $\nabla_vu$, $G$ and $F$ belong to
$L^2_{\rm loc}((0,T)\times\R^6)$.

\step{Step 6: the energy inequality.}
Let
\[
E(t):=\frac12\iint_{\R^6}u_+(t,x,v)^2\dd x\dd v .
\]
By \eqref{usupp} we apply Lemma~\ref{lem:pos}. Since $\nabla_vu_+=\mathbf 1_{\{u>0\}}\nabla_vu$, we get for a.e.\ $t\in(0,T)$
\[
E(t)=\int_0^t\!\!\iint_{\R^6}\Big(-\ba^g\nabla_vu_+\cdot\nabla_vu_++\beta^g\cdot\nabla_v\frac{u_+^2}{2}+F\,u_+\Big)\dd x\dd v\dd s .
\]
The first and the third term of the integrand are nonpositive, because $\ba^g\ge0$ and $F\le0$. For the second term we use that, for
a.e.\ $(s,x)$, the function $\zeta:=u_+^2/2$ belongs to $W^{1,1}(\R^3)$ and has compact support, so that by \eqref{div}
\[
\int_{\R^3}\beta^g\cdot\nabla_v\zeta\dd v=\int_{\R^3}\bc^g\,\zeta\dd v .
\]
Therefore, we have
\[
E(t)\le\frac12\int_0^t\!\!\iint_{\R^6}\bc^g\,u_+^2\dd x\dd v\dd s\le K'\int_0^tE(s)\dd s\qquad\text{for a.e.\ }t\in(0,T).
\]

\step{Step 7: conclusion.}
Since $E\in L^1(0,T)$, Gronwall's inequality gives $E=0$ a.e. Thus $u\le0$ a.e., and hence everywhere by continuity; that is,
\[
g\le\bar\phi+\varepsilon\psi\qquad\text{on }[0,T]\times\R^6 .
\]
Letting $\varepsilon\to0$, we obtain $g\le Me^{\beta_0T}\jv^{-k}$. This concludes the proof.
\end{proof}

\begin{proof}[Proof of Corollary~\ref{cor}]
Let $f_{\rm in}$, $f$ and $g$ be as in Corollary~\ref{cor}. We verify the hypotheses {\rm(A1)} and {\rm(A2)} of
Theorem~\ref{thm:main} for the pair $(f,g)$.

\step{Step 1: the initial datum.}
Since $g$ is uniformly continuous on $[0,T]\times\R^6$ and $g(t)\to f_{\rm in}$ pointwise as $t\searrow0$, we have
$g(0,\cdot,\cdot)=f_{\rm in}$, and $f_{\rm in}$ is uniformly continuous on $\R^6$.

\step{Step 2: the solution $f$.}
We use four facts from \cite{HST}.
\begin{itemize}
\item By \cite[Theorem 1.2]{HST}, $0\le f\in L^{\infty,k}$, and $f\in C^{2,\alpha}_{\rm kin,loc}((0,T]\times\R^6)$ solves \eqref{L}
pointwise. Hence $f$ is a weak solution, and therefore a distributional solution; see Section~\ref{sec:equation}.
\item By \cite[Proposition 3.1]{HST}, $f(t)\to f_{\rm in}$ uniformly as $t\searrow0$, because $f_{\rm in}$ is uniformly continuous.
\item By the weighted interior H\"older bound \cite[(2.17)]{HST} and the decay $|f|\le K\jv^{-k}$, the function $f$ is uniformly
continuous on $[t_0,T]\times\R^6$ for every $t_0>0$.
\item By \cite[Lemma 2.5]{HST} and the remark following it, which use the condition \cite[(1.7)]{HW}, the lower bound {\rm(A2)}
holds on $[0,T]$; see also \cite[(3.3)]{HW}.
\end{itemize}
The second and the third fact together show that $f$, with $f(0)=f_{\rm in}$, is uniformly continuous on $[0,T]\times\R^6$.

\step{Step 3: the solution $g$.}
By assumption, $g\ge0$ is a uniformly continuous weak solution in $L^{\infty,\ell}$ with $\ell=5+\gamma+\eta$. In particular, $g$ is
a distributional solution; see Section~\ref{sec:dual}. Since $g(0,\cdot,\cdot)=f_{\rm in}\in L^{\infty,k}(\R^6)$,
Lemma~\ref{lem:decay} gives
\[
g\in L^{\infty,k}([0,T]\times\R^6) .
\]

\step{Step 4: conclusion.}
By Steps 1--3, the functions $f$ and $g$ satisfy {\rm(A1)--(A2)} with $K:=\|f\|_{L^{\infty,k}}+\|g\|_{L^{\infty,k}}$. Hence
Theorem~\ref{thm:main} gives $f=g$ on $[0,T]\times\R^6$.
\end{proof}

\section*{Acknowledgments}
This work was partially supported by \rev{the Simons Foundation Travel Support for Mathematicians (No.~0007730) and by the AI for Science Award from Anthropic.} The author acknowledges the use of AI tools.

\end{document}